\documentclass[12pt,reqno]{amsart}

\numberwithin{equation}{section}

\usepackage[final, 
color]{showkeys}
\definecolor{refkey}{gray}{.85}
\definecolor{labelkey}{gray}{.85}

\usepackage{AKstyle}

\usepackage{caption}
\usepackage[labelfont=rm]{subcaption}

\usepackage[pagebackref=true, colorlinks]{hyperref}

\hypersetup{pdffitwindow=true,linkcolor=blue,citecolor=blue,urlcolor=blue,menucolor=blue}

\usepackage{comment}

\newcommand\faS{Y}
\newcommand\xiS{X}

\makeatletter
\let\orgdescriptionlabel\descriptionlabel
\renewcommand*{\descriptionlabel}[1]{%
  \let\orglabel\label
  \let\label\@gobble
  \phantomsection
  \edef\@currentlabel{#1}%
  \let\label\orglabel
  \orgdescriptionlabel{#1}%
}
\makeatother

\begin{document}

\author{Jean Bourgain}
\thanks{JB is partially supported by NSF grant DMS-0808042.}
\email{bourgain@ias.edu}
\address{School of Mathematics, IAS, Princeton, NJ}

\author{Alex Kontorovich}
\thanks{AK is partially supported by
an NSF CAREER grant DMS-1254788, an Alfred P. Sloan Research Fellowship, a Yale Junior Faculty Fellowship, and support at IAS from The Fund for Math and The Simonyi Fund.}
\email{alex.kontorovich@yale.edu}
\address{Department of Mathematics, Yale University, New Haven, CT}
\curraddr{Department of Mathematics, Rutgers University, New Brunswick, NJ}

\title[Low-Lying Fundamental Geodesics]
{Beyond Expansion II: Low-Lying Fundamental 
Geodesics}

\begin{abstract}
A closed geodesic on the modular surface is ``low-lying'' if it does not travel ``high'' into the cusp. It is ``fundamental'' if it corresponds to an element
in the class group of a real quadratic field. 
We prove the existence of infinitely many low-lying fundamental geodesics, answering a question of Einsiedler-Lindenstrauss-Michel-Venkatesh. 
\end{abstract}
\date{\today}
\subjclass[2010]{}
\maketitle
\tableofcontents


\section{Introduction}\label{sec:intro}

In this paper, we
answer a question 
of Einsiedler-Lindenstrauss-Michel-Venkatesh
on the abundance  of ``low-lying'' closed geodesics on the modular surface which are ``fundamental''; see the definitions below.
The main difficulty is to produce  a strong ``level of distribution'' for a particular set coming from a ``thin orbit.''

\subsection{Statement of the Main Theorem}\

Let $D>0$ be a fundamental discriminant, that is, the discriminant of a real quadratic field $K_{D}=\Q(\sqrt D)$, and let $\sC_{D}$ be the class group of $K_{D}$, with class number $h_{D}=|\sC_{D}|$. To each class $\g\in\sC_{D}$, we associate in the standard way a closed geodesic (by abuse of notation also called $\g$) in the unit tangent bundle 
$$
\cX \ :=\ \PSL_{2}(\Z)\bk\PSL_{2}(\R) \ \cong \  T^{1}(\PSL_{2}(\Z)\bk\bH) 
$$ 
of the modular surface. 
Not every closed geodesic 
on $\cX$
corresponds to an element of the class group of a real quadratic field; we call those that do {\bf fundamental}. The following rank-one question
arose around 2004 in the work of Einsiedler-Lindenstrauss-Michel-Venkatesh on higher rank analogues of Duke's Theorem (see \cite[\S1.5]{ELMV2009} and the discussion below).

\begin{question}\label{q:1}
Does there exist a compact subset $\cY\subset \cX$ which contains infinitely many fundamental geodesics?
\end{question}

Geodesics confined  to a compact region obviously never enter ``high'' in the cusp, and hence cannot equidistribute in $\cX$; we refer to 
these
as {\bf low-lying} (in $\cY$).
A natural set of candidate such, as observed by Sarnak, are the geodesics coming from Markov triples   (see \cite[pp. 226, 234]{Sarnak2007}), the difficulty being to understand when these are fundamental. 
(In fact, this very question initiated the study of the Affine Sieve 
\cite{BourgainGamburdSarnak2006, BourgainGamburdSarnak2010, SalehiSarnak2011
}.)
While we are unable to show the infinitude of fundamental Markov geodesics (which, if they exist, are extremely rare \cite{Zagier1982}), our
main
 goal 
 (see \thmref{thm:1} below)
 is to give an affirmative
 answer to
 \qref{q:1},
 in a strong quantitative sense.

Before stating our result, we
 put
\qref{q:1}
in perspective, by first recalling Duke's equidistribution theorem. 
Let $\mu_{\cX}$ be the probability Haar measure on $\cX$, and associate to each class $\g\in\sC_{D}$ (or rather, the corresponding geodesic) the probability arc-length measure $\mu_{\g}$.
Then
Duke's theorem \cite{Duke1988} asserts the equidistribution of $\mu_{\g}$'s to $\mu_{\cX}$ {\it on average} over $\sC_{D}$, for large discriminant:
\be\label{eq:Duke}
\frac1{h_{D}}\sum_{\g\in\sC_{D}}\mu_{\g}\ 
\xrightarrow{\ \text{weak } * \ }
 \ \mu_{\cX}
,\qquad\qquad\text{as $D\to\infty$.}
\ee
The goal of asking \qref{q:1} is to try to understand to what extent is it necessary to average over $\sC_{D}$ in \eqref{eq:Duke}, or whether perhaps 
the equidistribution already happens on the level of individual closed geodesic orbits (as is expected in higher rank analogues from rigidity phenomena conjectured by Cassels/Swinnerton-Dyer, Furstenberg, Margulis, etc). This question turns out to be quite subtle, as we indicate below. 

 \begin{figure}
        \begin{subfigure}[t]{.43\textwidth}
\includegraphics[width=1.95in]{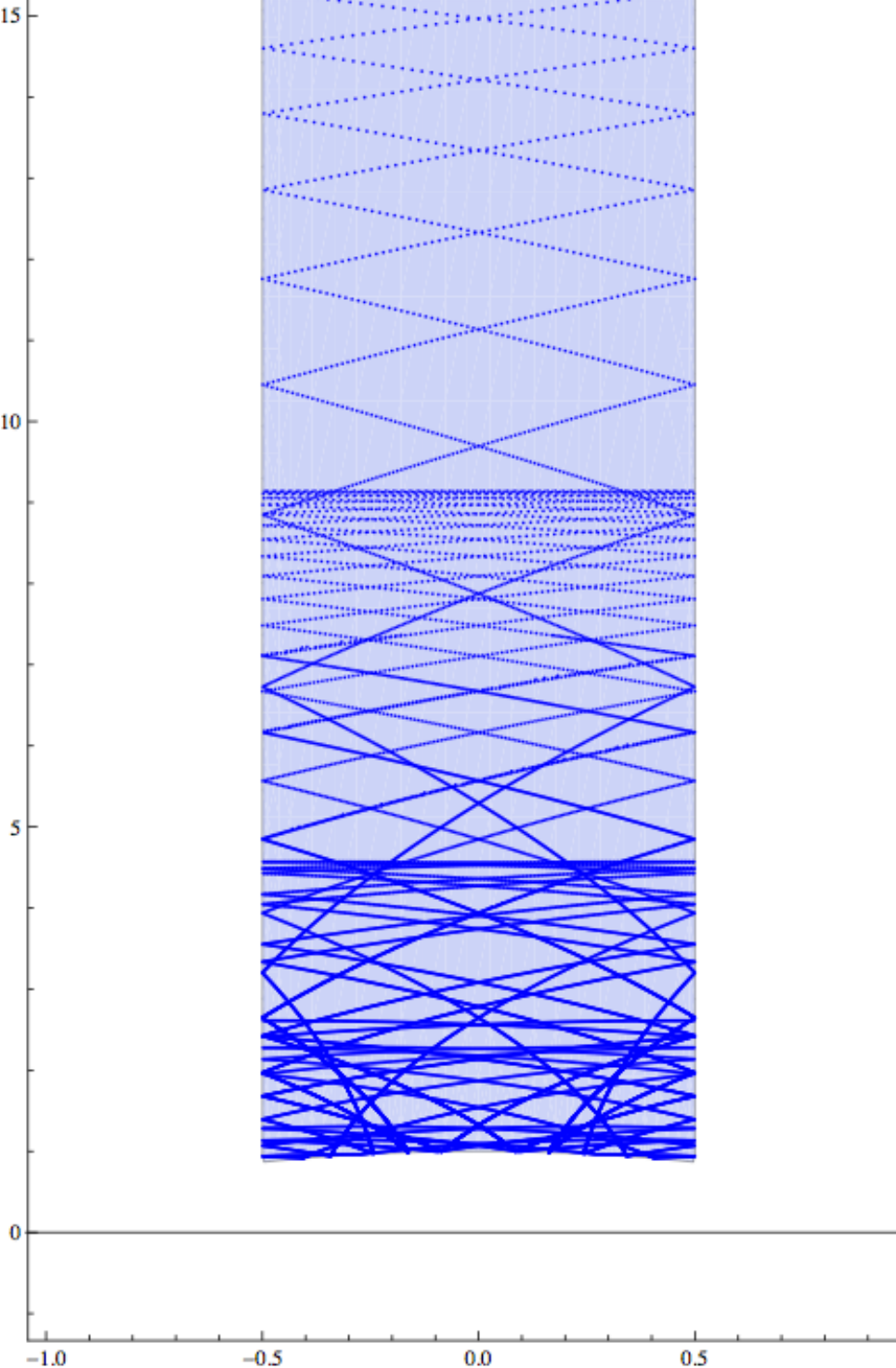}
                \caption{$D=1337$, $h_{D}=1$\\
                \phantom{nn} \protect\includegraphics[width=.1in]{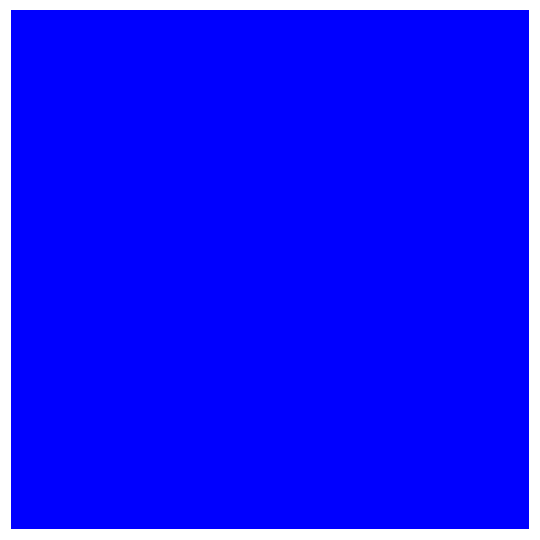} $=\{[19,27,-8]\}$}
                \label{fig:1a}
        \end{subfigure}%
\qquad
        \begin{subfigure}[t]{.5\textwidth}
                \includegraphics[width=1.95in]{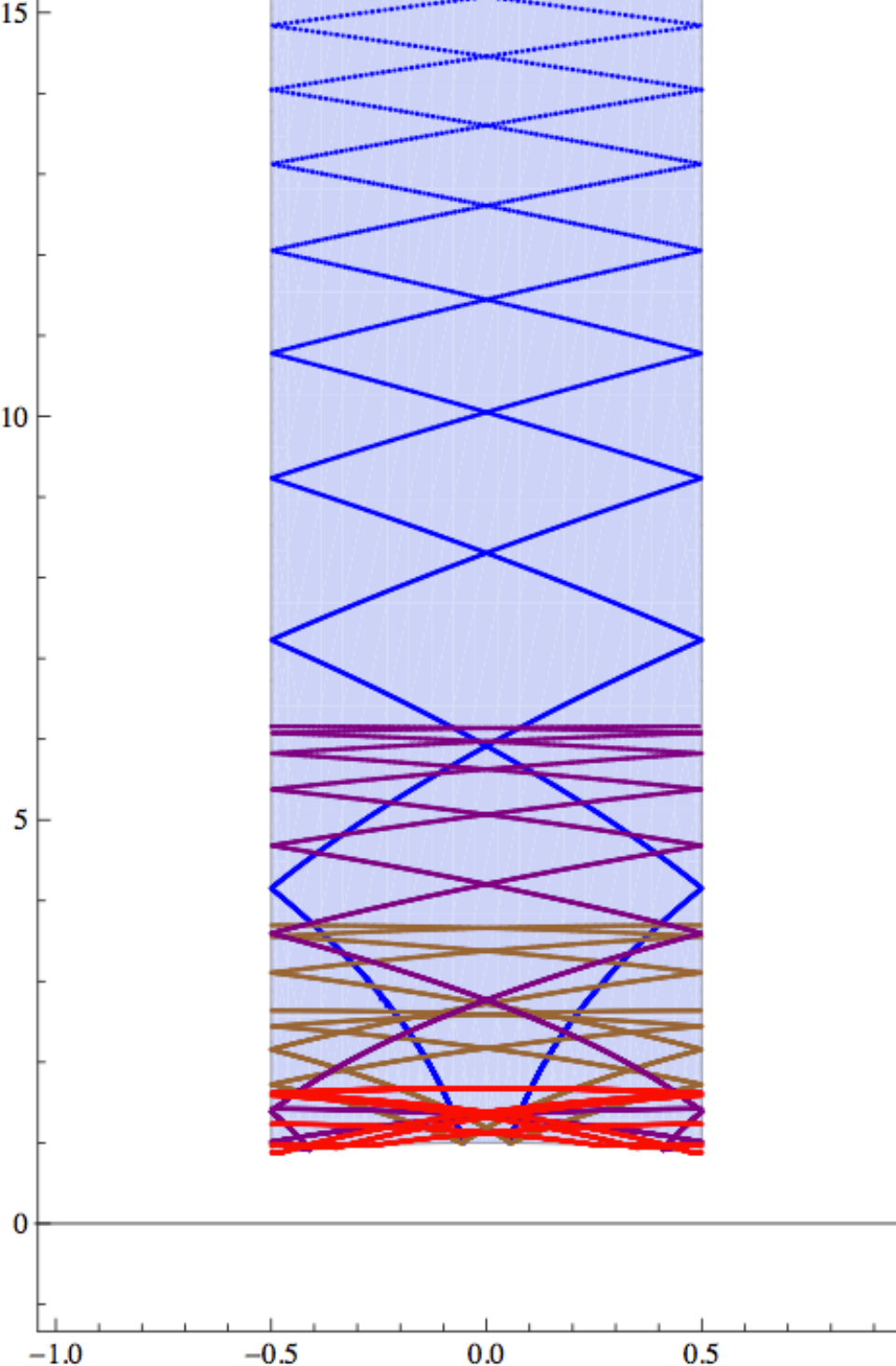}
                \caption{$D=1365$, $h_{D}=4$\\
                \protect\includegraphics[width=.1in]{SqBlue.pdf} $=\{[35,35,-1]\}$,
                \protect\includegraphics[width=.1in]{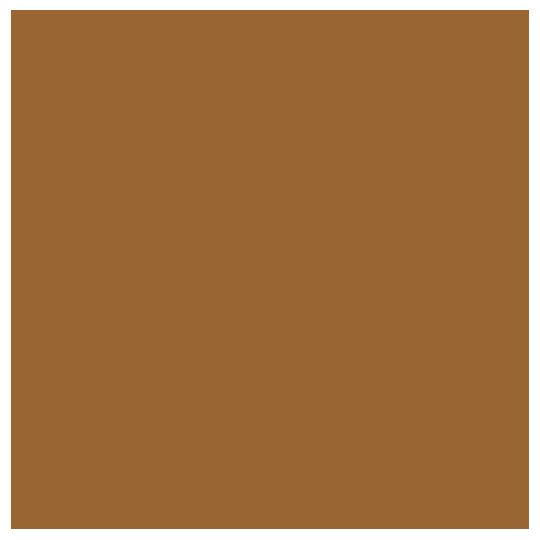} $=\{[7,35,-5]\}$,\\
                \protect\includegraphics[width=.1in]{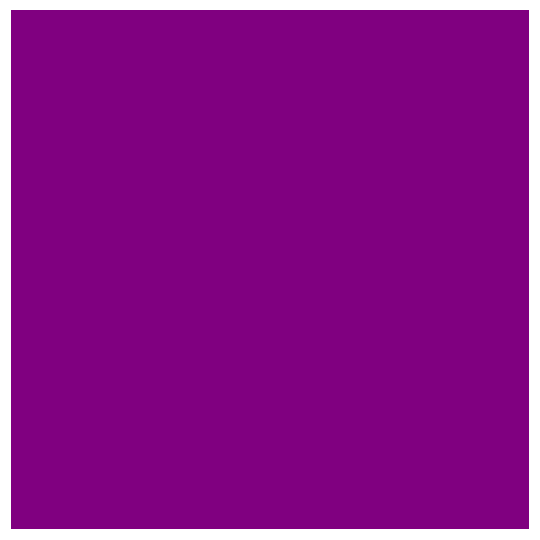} $=\{[23,33,-3]\}$,
                \protect\includegraphics[width=.1in]{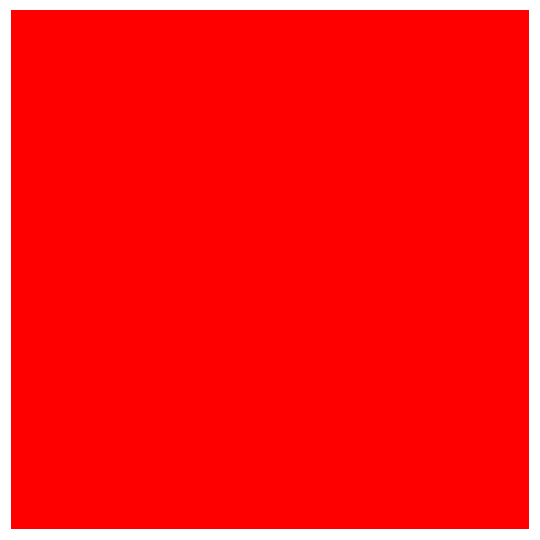} $=\{[19,23,-11]\}$
                }
                \label{fig:1b}
        \end{subfigure}

\caption{Fundamental geodesics in $\sC_{D}$.}
\label{fig:1}
\end{figure}

It will be instructive for the sequel to keep in mind the two examples illustrated in \figref{fig:1}. First recall some basic notions. Write $[A,B,C]$ for the binary quadratic form $Ax^{2}+Bxy+Cy^{2}$, and $\{[A,B,C]\}$ for the corresponding
 class. 
The discriminant 
\be\label{eq:discIs}
D\ = \ B^{2}-4AC,
\ee 
assumed throughout to be positive, is fundamental if
  either 
$D$  is square-free (in which case $D\equiv1\mod4$), or $D\equiv0\mod4$, 
in which case
$D/4$ is square-free and
$D/4\equiv2,3(\mod 4)$. 
%
%
 To associate a closed geodesic to a class $\g=\{[A,B,C]\}$, connect the two
  real Galois-conjugate
 roots 
\be\label{eq:geoGaDef}
\ga, \ \bar\ga\
=\
{-B\pm\sqrt D\over 2A},
\ee
by a geodesic in $T^{1}\bH$, and project modulo $\PSL_{2}(\Z)$. 
 We first consider the situation in

\phantomsection
{\bf Example I: $D=1337=7\times191$, \figref{fig:1a}.}\label{ex:1} While infinitely many closed geodesics are defined over $K=\Q(\sqrt{1337})$, only one of them is fundamental, because $K$ has class number one, $h_{1337}=1$. This one fundamental geodesic, corresponding to the class $\{[19,27,-8]\}$, is not particularly ``low-lying'' (of course that depends on one's choice of the compact region $\cY$), illustrating the difficulty of \qref{q:1}. 

In fact,
whenever
the field $K_{D}$ has class number one
 (as conjecturally happens  infinitely often),
there is
obviously
 no averaging in
Duke's theorem
 \eqref{eq:Duke}, and the one geodesic in the class is individually becoming equidistributed. 
Moreover, Popa's refinement \cite{Popa2006} of Waldspurger's theorem, together with a subconvex bound for certain Rankin-Selberg $L$-functions (see \cite{HarcosMichel2006}), implies that, as long as the class number is not too large, 
$$
h_{D} \ < \ D^{\eta}
$$ 
for some small $\eta>0$, then {\it every} geodesic in $\sC_{D}$ is {\it individually} becoming equidistributed; that is, without averaging over $\sC_{D}$ as in \eqref{eq:Duke}. Assuming the Lindel\"of hypothesis (which is a consequence of GRH) for such $L$-functions,  Popa's work implies that the exponent $\eta$ can be taken as large as $1/4-\vep$. Meanwhile, it is widely believed  that the same individual equidistribution holds with $\eta$ as large as $1/2-\vep$, for any fixed $\vep>0$.
So to even have a chance of seeing any non-equidistributing behavior (as in \thmref{thm:1}), one must take the class number almost as large as possible, 
\be\label{eq:largeH}
h_{D}\ > \ D^{1/2-o(1)}.
\ee

On the other hand, such discriminants should be quite rare. Indeed, it is 
a long-standing open problem 
 that the average class number satisfies, crudely,
$$
\sum_{0<D<T\atop \text{fundamental}}h_{D}\ 
\overset{?}=  \ T^{1+o(1)},\qquad\qquad{(T\to\infty)}.
$$
If true,
 this  estimate
 and the above heuristic
  would imply that there can only be very few discriminants with such large class number,
 \be\label{eq:fewDs}
\#\left\{0<D \text{ fundamental} <T\ : \ h_{D}>D^{1/2-o(1)}\right\}\ \overset{?}< \ T^{1/2+o(1)}.
 \ee 

Despite this rarity,
there
does exist a standard way of making large class numbers, namely, by considering discriminants of the special form 
\be\label{eq:DtoT}
D\ = \ t^{2}-4
.
\ee 
%
%
%
%
%
Then the fundamental solution to the Pellian equation
$$
T^{2}-DS^{2}\ = \ 4
$$
is $(T,S)=(t,1)$, whence the fundamental unit $\gep_{D}$ is as small as possible,
$$
\gep_{D}\ = \ {t+\sqrt D\over 2}
\asymp
\sqrt D.
$$
Dirichlet's Class Number Formula and Siegel's (ineffective) Theorem then give, crudely, that
$$
h_{D}\ = \ \sqrt D\ { L(1,\chi_{D})\over \log \gep_{D}}\ > \ D^{1/2-o(1)}
.
$$
Not surprisingly, we will be looking for low-lying (and hence non-equidistributing) behavior among fundamental discriminants of the special form \eqref{eq:DtoT}. 
(And since there are about $\sqrt T$ such up to $T$, we
 are not losing too much from \eqref{eq:fewDs}.) 
This brings us to:

\phantomsection\label{ex:2}
{\bf Example II: $D=1365=37^{2}-4= 3\times5\times7\times13$, \figref{fig:1b}.} 
The class number is $h_{1365}=4$, and the behavior of the four
fundamental geodesics defined over $\Q(\sqrt{1365})$ 
varies dramatically. The identity element of the class group $\sC_{1365}$
is 
the class $\{[35,35,-1]\}$, 
and the corresponding geodesic
shoots high up into the cusp;
meanwhile the geodesic corresponding to $\{[19,23,-11]\}$ is very low-lying, not reaching above $\Im z=2$. 
Nevertheless, 
the four geodesics
taken {\it together} 
equidistribute 
about
as well as the one geodesic in \hyperref[ex:1]{Example I}, beautifully illustrating why one must {\it  average} over $\sC_{D}$ for Duke's theorem \eqref{eq:Duke} to hold. 
\\

Returning to \qref{q:1}, we may now state our main result.


\begin{thm}\label{thm:1}
There exist infinitely many low-lying fundamental geodesics. More precisely, for each 
$\gep>0$, there is a compact region $\cY=\cY(\gep)\subset\cX$, and a set $\sD=\sD(\gep)$ of positive fundamental discriminants, such that:
\begin{enumerate}
\item
for each $D\in\sD$, many of the geodesics in the corresponding class group are low-lying:
\be\label{eq:sCsize}
\#\{\g\in\sC_{D}:\g\subset\cY
\}\ > \
 |\sC_{D}|^{1-\gep}
 ,
\ee
and moreover, 
\item
compared to \eqref{eq:fewDs},
there are many discriminants in $\sD$:
\be\label{eq:sDsize}
\#\sD\cap[1,T]
\ >
\ T^{1/2-\gep},\qquad\quad T\to\infty.
\ee
\end{enumerate}
\end{thm}

There are (at least) two ways to interpret this result. One can let $\gep\to0$, so that the inequalities \eqref{eq:sCsize}--\eqref{eq:sDsize} give
more and more ``low-lying'' 
fundamental
geodesics; 
unfortunately the compact region
$\cY(\gep)$ will then approach $\cX$,
giving less and less meaning to ``low-lying.''
Alternatively, one can let $\gep$ be a fixed  constant, 
say, $\gep=1/100$; then $\cY$ is a fixed region 
containing infinitely many fundamental geodesics, giving an affirmative answer to \qref{q:1}.

Again, in light of \eqref{eq:fewDs},
the estimate
\eqref{eq:sDsize} is almost sharp.
In \appref{sec:C}, we show (by more-or-less standard ergodic theoretic techniques, combining Duke's theorem and mixing) that \eqref{eq:sCsize} is also essentially sharp, in the following sense.
\begin{thm}\label{thm:sharp}
For any 
compact region $\cY\subset\cX$, there is an $\gep=\gep(\cY)>0$ so that
$$
\#\{\g\in\sC_{D}:\g
\subset\cY
\}\ < \
|\sC_{D}|^{1-\gep}
 ,
$$
as $D\to\infty$ through all non-square integers. 
\end{thm}


\subsection{Ingredients}\

We now describe some of the tools going into the proof of \thmref{thm:1}, beginning with a series of reformulations. 

\subsubsection{Step 1: Convert to Continued Fractions}\label{sec:ctdf}\

By the well-known connection  \cite{Humbert1916, Artin1924, Series1985} between continued fractions and the cutting sequence of the geodesic flow on $\cX$, 
the  condition that
a geodesic
 $\g$ be  low-lying can be 
reformulated as a
 Diophantine property on the
corresponding
visual   point 
$\ga$ in \eqref{eq:geoGaDef},
as follows.
Given
any $\ga\in\R$, write 
its
continued fraction expansion  as
$$
\ga=
a_{0}+\cfrac{1}{a_{1}+\cfrac{1}{a_{2}+\ddots}}
=
[a_{0},a_{1},a_{2},\dots]
,
$$
where $a_{0}\in\Z$ and $a_{j}\in\Z_{\ge1}$ for $j\ge1$;
the numbers $a_{j}$ are called
the
 {partial quotients} of $\ga$.
When $\ga$ is 
a
visual point of a closed geodesic, it is a quadratic irrational, and hence has an eventually periodic continued fraction expansion.
By applying a $\PSL_{2}(\Z)$ action, we may assume that $\ga$ is reduced, meaning that $-1<\bar\ga<0<1<\ga$;
then the continued fraction expansion of $\ga$ is exactly (as opposed to eventually) periodic.

Cutting off the cusp of $\cX$ at some height $\cC<\infty$  leaves a compact region $\cY=\cX\cap\{\Im(z)\le \cC\}$, and 
the condition that a geodesic $\g$ is low-lying (in $\cY$) is essentially equivalent to its visual point $\ga$
having all partial quotients bounded by some $\cA=\cA(\cC)<\infty$.

To illustrate this fact, we return for a moment to \figref{fig:1}. In \hyperref[ex:1]{Example I}, 
the one fundamental  geodesic in the class $\sC_{1337}$   corresponds to a reduced visual point $\ga$ having the continued fraction expansion:
$$
\{[19,27,-8]\}\quad\leftrightsquigarrow\quad[\overline{1, 1, 2, 17, 1, 8, 5, 8, 1, 17, 2, 1, 1, 3, 1, 35, 1, 3}].
$$
The large partial quotient $35$ is responsible for the high excursion of the geodesic in \figref{fig:1a}.

Meanwhile, the four geodesics in \hyperref[ex:2]{Example II} correspond to the continued fraction expansions:
\beann
\{[35,35,-1]\}&\leftrightsquigarrow& [\overline{1, 35}],
\\
\{[7,35,-5]\}&\leftrightsquigarrow& [\overline{5, 7}],
\\
\{[23,33,-3]\}&\leftrightsquigarrow& [\overline{1, 1, 1, 11}],
\\
\{[19,23,-11]\}&\leftrightsquigarrow& [\overline{1, 1, 1, 2, 1, 2}],
\eeann
and the very small partial quotients of the last of these explains the corresponding low-lying geodesic in \figref{fig:1b}.

To ensure that a fundamental geodesic is low-lying, one can try to force its visual point to 
have all partial quotients bounded by some height $\cA<\infty$. Alternatively, 
one can first consider all reduced quadratic irrationals with partial quotients bounded by  $\cA$, and try to understand 
when these
come from fundamental geodesics. We will take the latter approach, which turns out to be a sieving problem on a certain ``thin orbit.''

\subsubsection{Step 2: Convert to Thin Orbits}\

It is elementary that
the matrix
\be\label{eq:Majs}
\g
=
\mattwo{a_{0}}110
\mattwo{a_{1}}110
\cdots
\mattwo{a_{\ell}}110
\ee
fixes the quadratic irrational
\be\label{eq:gaDecomp}
\ga
=[\overline{a_{0},a_{1},\dots,a_{\ell}}]
,
\ee
thus
 converting
questions on continued fractions into ones about matrix products. In particular, we will be interested in traces of matrices of the form \eqref{eq:Majs}, in light of the following

\begin{lem}\label{lem:ELMVequiv}
A sufficient condition for a closed geodesic $[\g]$ to be fundamental is that 
\be\label{eq:ELMVequiv}
\tr(\g)^{2}-4\text{ is square-free.}
\ee
\end{lem}

Note that the corresponding discriminant $D=\tr(\g)^{2}-4$ is then fundamental and of the special form \eqref{eq:DtoT}.
Here is a quick proof:
The fixed points of $\g=\mattwos abcd$ are easily seen to be
\be\label{eq:gaToTr}
\ga,\ \bar\ga \ = \ {a-d\pm\sqrt{\tr(\g)^{2}-4}\over 2c}.
\ee
Now assume that $D:=\tr(\g)^{2}-4$ is square-free.
Comparing \eqref{eq:gaToTr} to 
\eqref{eq:geoGaDef}, we set $B:= d-a$ and $A:=c$. Solving \eqref{eq:discIs} for $C$ gives $C=-b$. Then the equivalence class of the form $Q=[A,B,C]$ has fundamental discriminant $D$, and hence corresponds to the fundamental geodesic $[\g]$, as desired.

 Thus to study  the traces of  matrix products of the form \eqref{eq:Majs} with 
 all
 $a_{j}\le \cA$, we should introduce the 
semigroup 
of finite products of such matrices,\footnote{The superscript $+$ in \eqref{eq:cGcAdef} denotes generation as a semigroup, that is, no inverses.}
\be\label{eq:cGcAdef}
\cG_{\cA}:=
\<
\mattwo a110
\ 
:
\ 
a\le \cA
\>^{+}
\quad\subset\quad\GL_{2}(\Z)
.
\ee
Preferring to work in $\SL_{2}$, we immediately pass to the even-length (determinant-one) sub-semi-group
\be\label{eq:GcAis}
\G_{\cA}
\ :=\ 
\cG_{\cA}\cap\SL_{2}(\Z)
,
\ee
which is (finitely) generated by the products $\mattwos a110\cdot\mattwos b110$, for $a,b\le \cA$.

The reason we call $\G_{\cA}$ ``thin'' is the following. 
Let $N$ be a growing parameter, and let
\be\label{eq:BNdef}
B_{N}\subset\SL_{2}(\R)
\ee 
be a ball about the origin of size $N$ in
the Frobenius norm 
$$
\|g\|^{2}=\tr(g{}^{t}g).
$$
A theorem of Hensley 
\cite{Hensley1989
} 
states
that
\be\label{eq:Hensley}
\#\{\G_{\cA}\cap B_{N}\}
\ \asymp_{\cA}   \ N^{2\gd_{\cA}},\qquad(N\to\infty),
\ee
where  $\gd_{\cA}$ is the Hausdorff dimension of the Cantor-like fractal of all numbers with partial quotients bounded by $\cA$: 
\be\label{eq:gdDef}
\gd_{\cA} \ := \ \text{H.dim}\,\{[a_{0},a_{1},a_{2},\dots]\ : \ a_{j}\le \cA\} \  \ \in(0,1)
.
\ee
This dimension has been studied extensively,
and it is known
\cite{Hensley1992} that  it can be made arbitrarily close to $1$ by taking $\cA$ large,
\be\label{eq:gdAlarge}
\gd_{\cA} \ = \ 1 - {6\over \pi^{2}\cA}+o\left(\frac1\cA\right), \qquad\qquad(\cA\to\infty).
\ee
On the other hand, 
the set of $\Z$-points in the Zariski closure of $\G_{\cA}$
is just  $\SL_{2}(\Z)$,
 and it is classical that $\#\{\SL_{2}(\Z)\cap B_{N}\}\asymp N^{2}$, instead of the much ``thinner'' count $N^{2\gd_{\cA}}$ as in \eqref{eq:Hensley}. 

In light of \lemref{lem:ELMVequiv} and \eqref{eq:Hensley},
the main
\thmref{thm:1} will follow without much effort from the following
\begin{thm}\label{thm:tr2m4}
Many elements in $\G_{\cA}$ have traces satisfying \eqref{eq:ELMVequiv}. More precisely,
for any $\eta>0$, there is an $\cA=\cA(\eta)<\infty$ such that
\be\label{eq:tr2m4}
\#\{
\g\in\G_{\cA}\cap B_{N}:\tr(\g)^{2}-4\text{ is square-free}\}
>
N^{2\gd_{\cA}-\eta}
,
\ee
as $N\to\infty$.
\end{thm}

The problem is thus reduced to

\subsubsection{Step 3:  Try to Execute a Sieve }\

This subsection is purely expository and heuristic; we
will give a 
  rough discussion of the sieving procedure, 
deferring
the precise (and somewhat technical) 
statements to \secref{sec:setup}. 

To sieve for square-free values of $\tr(\g)^{2}-4$, we need to understand their distribution modulo $\fq$, as $\g$ ranges roughly in $\G_{\cA}\cap B_{N}$, taking $\fq$ as large as possible relative to $N$. 
Since $\tr(\g)^{2}-4$ is of size $N^{2}$ when $\g$ is of size $N$, we introduce the parameter
\be\label{eq:Tdef}
T=N^{2}.
\ee
Letting $\gb(\fq)$ be the proportion of matrices in $\SL_{2}(\fq)$ for which the equation $\tr(\g)^{2}-4\equiv0(\mod \fq)$ holds,  one might expect that
$$
r_{\fq}(T) \ := \
\sum_{\g\in\G_{\cA}\cap B_{N}\atop\tr(\g)^{2}-4\equiv0(\fq)}1 \
- \
\gb(\fq)
\sum_{\g\in\G_{\cA}\cap B_{N}}1
$$
is a ``remainder'' term, which should be ``small'' in the following sense. 
We would like
that for some large $\cQ$, these remainders summed up to $\cQ$ still do not exceed the total size,
\be\label{eq:levRough}
\sum_{\fq<\cQ}|r_{\fq}(N)|
=
o
\left(
\#(\G_{\cA}\cap B_{N})
\right)
.
\ee
If this can be rigorously established, then we call $\cQ$
   a {\bf level of distribution} (for $\cA$).
Note that this is not a quantity intrinsic to our problem, but rather a function of what one can prove.
The larger this quantity, the more control one has on the distribution of the set of traces on such arithmetic progressions.
     If moreover \eqref{eq:levRough} can  be 
     confirmed
     with $\cQ$ as large as a power of the parameter $T$, 
\be\label{eq:cQNga}
\cQ=T^{\ga},\qquad\qquad
 \ga>0,
\ee 
then we say that $\ga$ is an {\bf exponent of distribution} for $\cA$.

The by-now ``standard'' Affine Sieve procedure applies in this context, and produces {\it some} (weak) exponent of distribution $\ga>0$. 
We briefly sketch the method now before explaining why it 
is insufficient
in our context.
A theorem of Bourgain-Gamburd-Sarnak \cite{BourgainGamburdSarnak2011} 
says
very roughly (see 
\thmref{thm:BGS2} for a precise statement) that
we do have equidistribution in $\G_{\cA}$ mod $\fq$, in the sense that there are constants 
\be\label{eq:gTgap}
\gT>0
\ee 
and $C<\infty$ so that, for all $\fq\ge1$ and all $\g_{0}\in\SL_{2}(\fq)$,
\be
\label{eq:BGSrough}
\left|
\sum_{\g\in\G_{\cA}\cap B_{N}\atop
\g\equiv\g_{0}(\fq)
}1 \
- \
{1\over |\SL_{2}(\fq)|}
\sum_{\g\in\G_{\cA}\cap B_{N}}1
%
%
\right|
\quad
\text{``} \ll\text{''}
\quad
\fq^{C}N^{2\gd-\gT}
,
\ee
where the implied constant is independent of $\g_{0}$ and $\fq$.
(We reiterate that the error in \eqref{eq:BGSrough} is heuristic only;  a statement of this strength is not currently known. That said, the true statement serves the same purpose in our application.) 
The positivity of $\gT$ in \eqref{eq:gTgap} {is} called the ``spectral gap'' or ``expander'' property of $\G_{\cA}$, and follows from a resonance-free region for the resolvent of a certain ``congruence'' transfer operator, see 
\secref{sec:BGS2}.
Summing \eqref{eq:BGSrough} over those $\g_{0}\in\SL_{2}(\fq)$ with $\tr(\g_{0})^{2}-4\equiv0(\fq)$, and then over $\fq$ up to $\cQ$, one proves \eqref{eq:levRough} with $\cQ=N^{\ga}$ and exponent of distribution
\be\label{eq:weakExp}
\ga  \ = \ \gT/C -\vep,
\ee
for any $\vep>0$.
(The value of $C$ may change from line to line.)

It turns out that this standard procedure is just shy of giving our main result! To successfully execute the sieve (that is, convert \eqref{eq:levRough} into \thmref{thm:tr2m4}), one needs the exponent of distribution 
to be strong enough to overcome the thinness of $\G_{\cA}$, in the sense that we need something like
\be\label{eq:gaToDel}
\ga \  > \ 10\cdot(1-\gd_{\cA});
\ee
see
\secref{sec:pftr2}.
The term on the right 
can be made arbitrarily small (cf. \eqref{eq:gdAlarge}), so it seems that by taking $\cA$ large enough, we should establish \eqref{eq:gaToDel}. Unfortunately, the spectral gap $\gT$ in \eqref{eq:gTgap} coming from the proof in \cite{BourgainGamburdSarnak2011} is {\it a priori} a function of $\cA$, 
and it is an extremely important open problem to understand its behavior with respect to $\cA$. Presumably $\gT$ should not deteriorate to zero as $\cA$ increases, but present methods are insufficient to show this, rendering the exponent \eqref{eq:weakExp} useless towards \eqref{eq:gaToDel}. 
Of course one can try to directly follow the proof in \cite{BourgainGamburdSarnak2011}, but then the $\cA$ dependence will be abysmal, and insufficient relative to \eqref{eq:gdAlarge} to produce the required inequality \eqref{eq:gaToDel}.

Rather than attacking this difficult problem head-on, 
we
circumvent
it as follows.

\subsubsection{Step 4: Prove an Exponent of Distribution Beyond Expansion}\

Instead of controlling the remainders \eqref{eq:levRough} 
using only expansion  \eqref{eq:BGSrough}, we seek to go beyond the direct procedure of the Affine Sieve, producing a stronger exponent of distribution 
to ensure that \eqref{eq:gaToDel} is satisfied.
We employ two novel techniques here, which appear in some form already in \cite{BourgainKontorovich2010, BourgainKontorovich2011, BourgainKontorovich2014,  BourgainKontorovich2014a, BourgainKontorovich2015a}. The first is to take inspiration from Vinogradov's method, 
replacing the full archimedean  ball $\G_{\cA}\cap B_{N}$ by a product of several such, which more readily captures the semigroup structure, and 
moreover allows development of techniques from estimating ``bilinear'' forms. 
The second innovation is, for larger values of $\fq$, to capture the condition $\tr(\g)^{2}-4\equiv0(\mod \fq)$ by {\it abelian} harmonic analysis, rather than the ``spectral'' method in \eqref{eq:BGSrough}. One then faces various exponential sums over our thin semigroup $\G_{\cA}$, but after some applications of Cauchy-Schwarz, one uses positivity to replace $\G_{\cA}$ by the full ambient group $\SL_{2}(\Z)$, allowing employment of more classical tools.
This loss is acceptable as long as the dimension $\gd_{\cA}$ of $\G_{\cA}$ is sufficiently near $1$, that is, as long as $\cA$ is large enough.
In the end, we are able to produce the strong exponent of distribution (for $\cA$ large) of
\be\label{eq:gaExp}
\ga=1/34.
\ee
In fact, our methods prove the exponent $\ga=1/32-\vep$ (and further refinements would give $\ga=1/8-\vep$), but \eqref{eq:gaExp}
 is already more than sufficient for \eqref{eq:gaToDel}; for ease of exposition, we will not strive for optimal exponents.

Inserting the strong exponent of distribution in \eqref{eq:gaExp} into \eqref{eq:gaToDel}, we are able to sieve down to square-free values of $\tr(\g)^{2}-4$, thus establishing \thmref{thm:tr2m4}, from which \thmref{thm:1} follows easily.

The proof of \thmref{thm:sharp} is by completely different methods, namely, a combination of mixing, Duke's theorem, and standard tools in ergodic theory.

\subsection{Organization}\

In \secref{sec:prelim}, we 
collect
some preliminaries needed in the sieve analysis. 
We spend \secref{sec:setup}  constructing the main ``bilinear'' set $\Pi\subset\G_{\cA}$ used for sieving, and stating the main sieving theorems. The main term is analyzed in  \secref{sec:main}, while the  errors are handled in \secref{sec:err}.
We prove the main sieving theorem (see \thmref{thm:sieve}) in \secref{sec:pfThm1}, and derive \thmref{thm:1} in \secref{sec:pfMain}.
Finally, the 
appendix proves \thmref{thm:sharp}.

\subsection{Notation}\

We use the following 
notation throughout. Set $e(x)=e^{2\pi i x}$ and $e_{q}(x)=e(\frac xq)$. We use the symbol $f\sim g$ to mean $f/g\to1$. The symbols $f\ll g$ and $f=O(g)$ 
are used 
interchangeably to mean the existence of an implied constant $C>0$ so that $f(x)\le C g(x)$ holds for all $x>C$; moreover $f\asymp g$ means $f\ll g\ll f$. 
 The letters $c$, $C$ denote  positive constants, not necessarily the same in each occurrence. 
Unless otherwise specified,  implied constants  may depend at most on $\cA$, 
which is treated as fixed. 
%
The letter $\vep>0$ is an arbitrarily small constant, not necessarily the same at each occurrence. When it appears in an inequality, the implied constant may also depend on $\vep$ without further specification.
The symbol $\bo_{\{\cdot\}}$ is the indicator function of the event $\{\cdot\}$. The trace of a matrix $\g$ is denoted $\tr\g$.
The number of divisors of $n$ is denoted $\tau(n)$. 
The greatest common divisor of $n$ and $m$ is written $(n,m)$ and their least common multiple is $[n,m]$.
The function $\nu(n)$ denotes the 
number of
distinct
 prime factors of $n$. 
The cardinality of a finite set $S$ is denoted $|S|$ or $\# S$.
The transpose of a matrix $g$ is written ${}^{t}g$. When there can be no confusion, we use the shorthand $r(q)$ for $r(\mod q)$. The prime symbol $'$ in $\underset{r(q)}{\gS} {}'$ means the range of $r(\mod q)$ is restricted to $(r,q)=1$. 

\subsection*{Acknowledgments}\

It is our pleasure to thank
  Tim Browning, Curt McMullen, Michael Rubinstein, Zeev Rudnick, and Peter Sarnak for illuminating conversations. The second-named author would like to thank the hospitality of the IAS, where much of this work was carried out.


\section{Preliminaries}\label{sec:prelim}

In this section, we state two results that are needed in the sequel, namely \propref{prop:aleph} and \propref{prop:SL2Bnd}. We recommend the reader skip the 
proofs on the first pass, instead proceeding directly to \secref{sec:setup}. 

\subsection{Expansion}\label{sec:BGS2}\

In this subsection, we make precise the ``expansion'' theorem heuristically stated in \eqref{eq:BGSrough}.
We will only require expansion for the fixed 
value
$\cA_{0}=2$, so as to make the expansion constants absolute, and not dependent on $\cA$; see 
the discussion after \eqref{eq:gaToDel} and
\rmkref{rmk:gTunif}.

To this end, let $\G_{2}\
%
\subset\SL_{2}(\Z)$
be the
  semigroup
  as  in \eqref{eq:GcAis}
  corresponding to $\cA_{0}:=2$. 
It is easy to see that $\G_{2}$ is free, that every non-identity matrix $\g\in\G_{2}$ is hyperbolic, and that
$$
\tr\g\ \asymp\ \|\g\|.
$$
Note also that the {\it group} $\<\G_{2}\>$ generated by the semigroup $\G_{2}$ is all of $\SL_{2}(\Z)$.
This
immediately implies that for 
any
$q\ge1$, 
the mod $q$ reduction of $\G_{2}$ is everything,
\be\label{eq:Gmodq}
\G_{2}(\mod q)\ \cong \ \SL_{2}(q)
.
\ee

 The following theorem is a consequence of the general expansion theorem proved by Bourgain-Gamburd-Sarnak in \cite{BourgainGamburdSarnak2011}.

\begin{thm}[\cite{BourgainGamburdSarnak2011}]\label{thm:BGS2}
Let  $
\G_{2}$ be the semigroup above. 
There exists
an absolute square-free integer
\be\label{eq:fBis}
\fB
\ge1
\ee
absolute constants $c,\ C
>0$, 
and an absolute ``spectral gap'' 
\be\label{eq:gTis}
\gT=\gT(\cA_{0})>0,
\ee
so that,
for any square-free $q\equiv0(\fB)$ and
any $\gw\in\SL_{2}(q)$, 
as $\faS\to\infty$,
we have
\bea
\label{eq:BGS2}
&&
\hskip-.5in
\#\left\{
\g\in\G_{2}\cap B_{\faS}
:
\g\equiv\gw(\mod q)
\right\}
\\
\nonumber
&=&
{|\SL_{2}(\fB)|
\over |\SL_{2}(q)|}
\bigg|
\left\{
\g\in\G_{2}\cap B_{\faS}
:
\g\equiv\gw(\mod\fB)
\right\}
\bigg|
\\
\nonumber
&&
\hskip1in
+
O\Big(
\#\left\{
\g\in\G_{2}
:
\|\g\|<\faS
\right\}
\cdot
\fE(\faS;q)
\Big)
,
\eea
where
\be\label{eq:fEbnd}
\fE(\faS;q
)
:=
\twocase{}
{e^{-c\sqrt{\log \faS}},}{if $q\le  C\log \faS,$}
{q^{C} \faS^{-\gT},}{if $q>C\log \faS.$}
\ee
 \end{thm}

\begin{rmk}
This theorem is 
proved 
in \cite[see Theorem 1.5]{BourgainGamburdSarnak2011} under the 
assumption
that, instead of
 $\G_{2}$, one is given a convex-cocompact 
sub{\it group} of $\SL_{2}(\Z)$. 
But the proof is the same when the group is replaced by our free semigroup $\G_{2}$;
we emphasize again that $\G_{2}$
  has no parabolic elements.
The error term \eqref{eq:fEbnd} is the consequence of a Tauberian argument applied to a resonance-free region \cite[see Theorem 9.1]{BourgainGamburdSarnak2011} of the form 
\be\label{eq:resolvent}
\gs>\gd_{\cA_{0}}-C\min\left\{1,{\log q\over \log (1+|t|)}\right\},\qquad 
\gs+it\in\C,
\ee
for the resolvent 
of a certain ``congruence'' transfer operator
, see \cite[\S12]{BourgainGamburdSarnak2011} for details. For small $q$, we only obtain a ``Prime Number Theorem''-quality error (given here in crude form), while for larger $q$, \eqref{eq:resolvent} is as good as a resonance-free strip. 
\end{rmk}

We have stated the result only for the 
case 
$\fB\mid q$. 
The distribution modulo $\fB$ cannot be obtained directly from present methods, 
even 
though
all reductions of $\G_{2}$ are surjective;  cf. \eqref{eq:Gmodq}. 
Nevertheless, one can construct a set which has the desired equidistribution for all $q$, as claimed in the following

\begin{prop}\label{prop:aleph}
Given any $\faS\gg1$, 
there is a non-empty subset 
$$
\aleph=\aleph(\faS)\ \subset \ \G_{2}\cap B_{\faS}
$$ 
so that,
for all square-free $q
$ and all $\fa_{0}\in\SL_{2}(q)$, 
\be\label{eq:alephED}
\left|
{\#\{\fa\in\aleph:\fa\equiv\fa_{0}(q)\}\over |\aleph|}-{1\over |\SL_{2}(q)|}\right|
\ll
\fE(\faS;q)
.
\ee
Here $\fE$ is given in \eqref{eq:fEbnd}.
\end{prop}

Note that we do not have particularly good control on the cardinality of $\aleph$;
regardless, the estimate \eqref{eq:alephED} is only 
nontrivial
if 
$$
q<\faS^{{\gT}/C}.
$$
The construction of the set $\aleph$ proceeds in a similar way to \cite[\S8]{BourgainKontorovich2014}; we  
give a sketch 
for the reader's convenience.

\pf[Proof (Sketch)]
Let the constants $\fB,$ $c,$ $C,$ and $\gT$ be as in 
\thmref{thm:BGS2}; they depend only on $\cA_{0}=2$, that is, they are absolute. 

Let $U$ be a parameter to be chosen later relative to $\faS$. Let 
$$
R:=|\SL_{2}(\fB)|\asymp1,
$$
and let
$$
\cS(U):=\{\g\in\G_{2}:\|\g\|<U\}
.
$$
From \eqref{eq:Hensley}, we have that
$$
\#\cS(U)\gg U^{2\gd_{2}},
$$
where $\gd_{2}=\gd_{\cA_{0}}$ is the corresponding Hausdorff dimension. 
Then by the pigeonhole principle, there exists some $\fs_{U}\in\cS(U)$ so that
$$
\cS'(U):=\{\g\in\G_{2}:\|\g\|<U,\ \g\equiv\fs_{U}(\mod\fB)\}
$$
has cardinality at least
$$
\#\cS'(U)\ge{1\over R}\#\cS(U)\gg U^{2\gd_{2}}.
$$

Observe that  the elements in
$$
\cS'(U)\cdot \fs_{U}^{R-1}
$$
are all congruent 
to
the identity mod $\fB$. 
Write $\SL_{2}(\fB)=\{\g_{1},\dots,\g_{R}\}$, and find $\fx_{1},\dots,\fx_{R}\in\G$ with
$$
\fx_{j}\equiv\g_{j}(\mod \fB).
$$
Such $\fx_{j}$ can be found of size
$$
\|\fx_{j}\|\ll1.
$$

For each $j=1,\dots,R$, let
$$
\cS_{j}'(U):=\cS'(U)\cdot \fs_{U}^{R-1}\cdot \fx_{j}
.
$$
This is a subset of $\G_{2}$, in which each element $\fs\in\cS_{j}'(U)$ has size as most 
$$
\|\fs\|
\ll
U^{R}
.
$$
Choose $U\asymp \faS^{1/R}$ so that all elements  $\fs\in\bigcup_{j}\cS_{j}'(U)$ 
satisfy
$$
\|\fs\|<\faS
.
$$
%
Then we claim that
$$
\aleph:=\bigsqcup_{j=1}^{R}\cS_{j}'(U)
$$
gives the desired special set.

Indeed, applying 
\thmref{thm:BGS2} gives that for each $j=1,\dots, R$, any $q
$ with $q\equiv0(\fB)$, and any $\gw\in\SL_{2}(q)$ with $\gw\equiv\fx_{j}(\fB)$, we have
\beann
&&
\hskip-.5in
\#\{\fs\in\cS_{j}'(U):\fs\equiv\gw(q)\}
\ =\
\#\{\fs\in\cS'(U):\fs\equiv\gw(\fs_{U}^{R-1}\fx_{j})^{-1}(q)\}
\\
&=&
\#\{\fs\in\cS(U):\fs\equiv\gw(\fs_{U}^{R-1}\fx_{j})^{-1}(q)\}
\\
&=&
{|\SL_{2}(\fB)|\over |\SL_{2}(q)|}
\#\{\fs\in\cS(U):\fs\equiv
\fs_{U}
(\fB)\}
+
O(
\fE(U;q)
)
\\
&=&
{|\SL_{2}(\fB)|\over |\SL_{2}(q)|}
\#\cS_{j}'(U)
+
O(
\fE(U;q)
)
.
\eeann

Then the sets $\cS_{j}'(U)$ each have good modular distribution properties for distinct residues mod $\fB$. Note that they also all have the same cardinality, namely that of $\cS'(U)$.
Moreover after renaming constants, we have that $\fE(U;q)\ll\fE(\faS;q)$.

The equidistribution \eqref{eq:alephED} is now clear for any $q\equiv0(\fB)$, while the same for other $q$
 is
 obtained by 
  summing over suitable arithmetic progressions. This completes the proof.
\epf


\subsection{An Exponential Sum over $\SL_{2}(\Z)$}\ 

In this subsection, we 
state an estimate, showing roughly 
 that there is cancellation in a certain 
exponential sum over $\SL_{2}(\Z)$ in a ball. 
We identify $\Z^{4}$ with $M_{2\times2}(\Z)$, and observe that for $A,B\in M_{2\times2}(\Z)$,
\be\label{eq:ZM2}
\tr(^{t}\hskip-2ptAB)=A\cdot B,
\ee
where the operation on the right is the dot product in $\Z^{4}$. 
We first give the following local result.

\begin{lem}\label{lem:modq}
For any square-free $q\ge1$ and any 
 vector $
\bs
\in
\Z^{4}$ 
with $(\bs,q)=1$,
we have
$$
\left|\sum_{\g\in\SL_{2}(q)}e_{q}(\g\cdot\bs)\right|
\ll
q^{3/2+\vep},
$$
for any $\vep>0$.
\end{lem}
\pf
We could appeal to Deligne, but in fact the estimate is elementary, involving only Weil's bound for Kloosterman sums.
The left hand side is multiplicative and $q$ is square-free, so we may consider just the case that $q=p$ is prime. Writing $\bs=(x,y,z,w)$, we may assume that, say,  $y\not\equiv0(p)$. Writing $\g\in\SL_{2}(p)$ as $\g=\mattwos abcd$, we break the sum according to whether or not $c\equiv0$. The former case contributes
$$
\sideset{}{'}\sum_{a(p)}
\sum_{b(p)}
e_{p}(ax+by+\bar a w)
=0,
$$
since $y\not\equiv0(p)$. When $c\not\equiv0$, we have:
\beann
&&
\hskip-1in
\sideset{}{'}\sum_{c(p)}
\sum_{a,d(p)}
e_{p}(ax+\bar c(ad-1)y+cz + d w)
\\
&=&
\sideset{}{'}\sum_{c(p)}
e_{p}(cz-\bar cy)
\sum_{a(p)}
e_{p}(ax)
\sum_{d(p)}
e_{p}(d(\bar c ay +  w))
\eeann
The $d$ sum vanishes except for the one value of $a\equiv-  c\bar yw(\mod p)$, in which case it contributes $p$. What remains is a Kloosterman sum in $c$, which is bounded by $2\sqrt p$, since $y\neq0(p)$. 
\epf

The next result
we record
is 
well-known, see, e.g.,
a special case of
 \cite[Theorem 2.9]{BourgainKontorovich2015a}.

\begin{lem}\label{lem:arch}
Let $\xiS\gg1$ be a growing parameter.
There exists a function 
$$
\vf_{X}:\SL_{2}(\R)\to\R\ge_{0}
$$ 
which approximates the indicator of an archimedean ball, 
by which we mean the following.
We have the lower bound
\be\label{eq:vfXlwr}
\vf_{X}(g) \ \ge \ 1
,
\ee
on 
$\|g\|\le X$,
and the upper bound
\be\label{eq:vfXupr}
\sum_{\xi\in\SL_{2}(\Z)}\vf_{X}(\xi) \  \ll \ X^{2}.
\ee
Moreover,
for any $q\ge1$, and any $\g_{0}\in\SL_{2}(q)$, 
\be\label{eq:vfXeq}
\sum_{\xi\in\SL_{2}(\Z)\atop\xi\equiv\g_{0}(q)}\vf_{X}(\xi)
 \ = \
\frac1{|\SL_{2}(q)|}
\sum_{\xi\in\SL_{2}(\Z)}\vf_{X}(\xi)
+
O(X^{3/2})
.
\ee
\end{lem}

The error term in \eqref{eq:vfXeq} comes from applying Selberg's $3/16$th spectral gap \cite{Selberg1965}; of course better estimates are now known \cite{KimSarnak2003}, but since we are not optimizing exponents, we will use the simplest results which suffice.

Combining the previous two lemmata, the main result of this subsection is the following

\begin{prop}\label{prop:SL2Bnd}
Let $\vf_{X}$ be as in \lemref{lem:arch}. For any square-free $q
\ge 1$, 
and
any 
 vector $
\bs
\in
\Z^{4}$ 
with $(\bs,q)=1$,
we have
for any $\vep>0$,
\be\label{eq:SL2Bnd}
\left|
\sum_{\xi\in\SL(2,\Z)}\vf_{\xiS}\left({\xi}\right)
e_{q}\left(
\xi\cdot \bs\right)
\right|
\ll
q
^{-3/2+\vep}\xiS^{2}
+
q^{3}
\,
\xiS^{3/2}
,
\ee
as $\xiS\to\infty$.
\end{prop}
\pf
Decompose the left side of \eqref{eq:SL2Bnd} according to the residue class of $\xi$ in $\SL_{2}(q)$, and apply \eqref{eq:vfXeq} and \eqref{eq:vfXupr}, giving
\beann
\left|
\sum_{\xi\in\SL(2,\Z)}\vf_{\xiS}\left({\xi}\right)
e_{q}\left(\xi\cdot \bs\right)
\right|
&=&
\left|
\sum_{\g\in\SL_{2}(q)}
e_{q}\left(\g\cdot \bs\right)
\sum_{\xi\in\SL(2,\Z)\atop\xi\equiv\g(q)}\vf_{\xiS}\left({\xi}\right)
\right|
\\
&\ll &
\left|
\sum_{\g\in\SL_{2}(q)}
e_{q}\left(\g\cdot \bs\right)
\right|
{X^{2}\over |\SL_{2}(q)|}
+O(q^{3}X^{3/2})
.
\eeann
The estimate follows from  \lemref{lem:modq}.
\epf


\section{%
Construction of $\Pi$ and 
the Sieving Theorem
}\label{sec:setup}

\subsection{Construction of the set $\Pi$}\

The first goal in this subsection is to construct an appropriate
subset $\Pi$  of $\G_{\cA}\cap B_{N}$
in
which to execute our sieve. Let
$\cA
<\infty
$ be 
fixed,
let
 $\G_{\cA}$ be the semigroup in \eqref{eq:GcAis},
 and let
$\gd_{\cA}$
be the
corresponding dimension \eqref{eq:gdDef},
assumed to be 
 near $1$. 
Since $\cA$ is fixed, we 
henceforth
drop the subscripts, writing $\G=\G_{\cA}$ and $\gd=\gd_{\cA}$. Recall also that implied constants may depend on $\cA$ without further specification.

Recall that $N$ is our main growing parameter,
 and let
\be\label{eq:Xis}
\xiS\ = \ N^{x},\quad  \faS\ = \ N^{y}, \quad Z\ = \ N^{z},\quad x,y,z\ > \ 0,
\ee 
be some parameters to be chosen later;
they will decompose $N$, in the sense that
\be\label{eq:NXYZ}
N\ =\ \xiS \faS Z, \qquad \quad\text{ that is, } \quad \quad x+y+z\ = \   1.
\ee
We think of $\xiS$ as  large, 
$\xiS>N^{1-\eta}$,
 and $\faS$ and $Z$ 
 as tiny.

Let $\aleph=\aleph(\faS)\subset\G_{2}\subset\G$ be the set constructed in 
\propref{prop:aleph}, and let 
\be\label{eq:gWdef}
\Xi_{0}\ :=\ \{\xi\in\G:\|\xi\|<\xiS\}, \qquad
\gW_{0}\ :=\ \{\gw\in\G:\|\gw\|<Z\}
\ee
be norm balls in $\G$. 
While we do not have good control on the size of $|\aleph|$,
recall from \eqref{eq:Hensley} that
\be\label{eq:gWsize1}
|\Xi_{0}|\ \asymp \ \xiS^{2\gd},\qquad\qquad
|\gW_{0}|\ \asymp\  Z^{2\gd}
.
\ee
We will want the products 
$$
\xi_{0}\cdot\fa\cdot\gw_{0}
$$ 
to be unique for $\xi_{0}\in\Xi_{0},$  $\fa\in\aleph,$  $\gw_{0}\in\gW_{0}$; since $\G_{\cA}$ is a free finitely-generated semigroup, this will be the case if the wordlength $\ell(\cdot)$ in the generators \eqref{eq:cGcAdef} is fixed in each norm ball. 
It is easy to see that the wordlength metric is commensurable to the log-norm,
\be\label{eq:ell}
\ell(\g)\ \asymp
\ \log\norm{\g}.
\ee
Then by the pigeonhole principle and \eqref{eq:gWsize1},  there is some $\ell_{X}\asymp\log X$ so that, defining
$$
\Xi \ 
:= \ \{\g\in\Xi_{0}:\ell(\g)=\ell_{X}\},
$$
we have
\be\label{eq:Xisize}
\#\Xi\ \gg\ {X^{2\gd}\over \log X}.
\ee
Similarly, there is some $\ell_{Z}\asymp \log Z$ so that, defining $\gW$ to be the subset of $\gW_{0}$ having wordlength exactly $\ell_{Z}$,
gives a set of cardinality
\be\label{eq:gWsize}
\#\gW\ \gg\ {Z^{2\gd}\over \log Z}.
\ee
Then the product 
\be\label{eq:PiDef}
\Pi \ := \ \Xi\cdot\aleph\cdot\gW
\ee
is a 
subset (and not a multi-set, since the products are unique) of $\G$. By \eqref{eq:NXYZ}, we clearly have
\be\label{eq:PiSub}
\Pi\ \subset \ \G\cap B_{100N}.
\ee
The set $\Pi$ will have our desired ``bilinear'' (in fact, multi-linear) structure, suitable for sieving.

\subsection{Statement of the Sieving Theorem}\ 
 
We can finally state the main sieving theorems for $\Pi$.

\begin{thm}\label{thm:AP}
Let $\Pi_{AP}$ denote the set of elements $\vp\in\Pi$ for which   $\tr(\vp)^{2}-4$ is ``almost prime,'' in particular having no ``small'' prime factors,
$$
\Pi_{AP}\ := \
\{
\vp\in\Pi\ : \ p\mid (\tr(\vp)^{2}-4)
\ \Longrightarrow\  p>
N^{1/350}\}
.
$$
Then
for
 any
 sufficiently small
  $\eta>0$,
 there is an
 $\cA=\cA(\eta)$ sufficiently large,
 and a choice of the parameters $X$, $Y$, $Z$ in \eqref{eq:Xis}, so that
  we have
\be\label{eq:PiAPbnd}
\#\Pi_{AP}
\ >\ 
N^{2\gd-\eta}
,
\ee
as $N\to\infty$.
\end{thm}

\thmref{thm:AP} will easily imply \thmref{thm:tr2m4}, and will itself be easily implied by the following ``level of distribution'' result.

Define the sifting sequence 
$$
\fA \ = \ \{a_{N}\},
$$ 
by
\be\label{eq:aNdef}
a_{N}(n)\ := \
\sum_{\vp\in\Pi}
\bo_{\{\tr^{2}(\vp)-4\ = \ n\}}
.
\ee
Note that, by \eqref{eq:PiSub}, $\fA$ 
has support
\be\label{eq:suppA}
\supp\fA\ \subset \ \{
n\ll T\}
,
\ee
where $T=N^{2}$, cf.  \eqref{eq:Tdef}. 
For a parameter $\cQ$ and any square-free $\fq<\cQ$, we define
\be\label{eq:fAfqDef}
|\fA_{\fq}|
\ := \
\sum_{n\equiv0(\fq)}a_{N}(n)
.
\ee
\begin{thm}\label{thm:sieve}
For any sufficiently small $\eta>0$, there is an $\cA=\cA(\eta)$ sufficiently large, so that the sequence $\fA$ has level of distribution 
\be\label{eq:ga110}
\cQ=T^{1/32-\eta
}.
\ee
More precisely, there is a multiplicative function $\gb:\N\to\R$
satisfying the
``quadratic sieve'' condition
\be\label{eq:gbEst}
\prod_{w\le p <z}
(1-\gb(p))^{-1}
\ \le \ 
C\cdot
\left({\log z\over \log w}\right)^{2}
,
\ee
for some $C>1$ and any $2\le w<z$;
and
 a decomposition
\be\label{eq:fAfq}
|\fA_{\fq}|\ =\ \gb(\fq)\cdot |\Pi|+r({\fq}),
\ee
so that, for all $K<\infty$,
\be\label{eq:rfqLevel}
\sum_{\fq<\cQ\atop\text{square-free}}|r(\fq)| \ \ll_{K} \ {|\Pi|\over \log^{K}N}.
\ee
Moreover, 
we can choose
\be\label{eq:Xchoice}
X=N^{1-\eta}
\ee 
in the decomposition \eqref{eq:PiDef} of
 $\Pi$, so that 
\be\label{eq:PiSize}
\#\Pi\ \gg \ N^{2\gd-\eta}.
\ee
\end{thm}
We now give a quick
\pf[Sketch of \thmref{thm:AP} assuming \thmref{thm:sieve}]\

The deduction is standard.
The content of the latter is that the sifting sequence $\fA$ has ``sieve dimension'' $\gk=2$, and any exponent of distribution $\ga<1/32$; this confirms  
the discussion below \eqref{eq:gaExp}.
Taking $\ga=1/34$, say (again, we are not striving for optimal exponents), and using the crudest Brun sieve, see, e.g. \cite[Theorem 6.9]{FriedlanderIwaniecBook}, one shows that
\be\label{eq:aNsieve}
\sum_{n\atop (n,P_{z})=1}a_{N}(n) \ \gg \  {|\Pi|\over (\log N)^{2}},
\ee
where $P_{z}=\prod_{p<z}p$ and $z$ does not exceed $T^{\ga/(9\gk+1)}=T^{1/646}=N^{1/323}$; we take $z=N^{1/350}$.
Of course any $n=\tr(\vp)^{2}-4$ coprime to $P_{z}$ has no prime factors below $z$. 
Then \eqref{eq:aNsieve} and \eqref{eq:PiSize}
confirm
 \eqref{eq:PiAPbnd} after renaming constants.
\epf
We focus henceforth on establishing \thmref{thm:sieve}. 

\subsection{The Decomposition }\

The decomposition
\eqref{eq:fAfq}
is determined as follows.
Inserting \eqref{eq:aNdef} into \eqref{eq:fAfqDef} gives
$$
|\fA_{\fq}|
\ =  \ 
\sum_{\vp\in\Pi}
\bo_{\{\tr^{2}(\vp)-4 \equiv  0(\fq)\}}
\ =  \ 
\sum_{\ft\mod\fq\atop\ft^{2}\equiv4(\fq)}
\sum_{\vp\in\Pi}
\bo_{\{\tr(\vp) \equiv  \ft(\fq)\}}
.
$$

Rather than applying expansion (that is, the analogue of \eqref{eq:BGSrough}) directly, we first capture the indicator function by {\it abelian} harmonic analysis, writing
$$
|\fA_{\fq}|
\ =  \ 
\sum_{\ft\mod\fq\atop\ft^{2}\equiv4(\fq)}
\sum_{\vp\in\Pi}
\frac1\fq
\sum_{q\mid\fq}\sideset{}{'}\sum_{r(q)}
e_{q}(r(\tr(\vp) -  \ft))
.
$$
Introducing a new parameter $Q_{0}<\cQ$, we obtain the decomposition \eqref{eq:fAfq}  from breaking the penultimate sum above 
according to whether $q<Q_{0}$ or not. To this end, we write
\be\label{eq:fAfqDecomp}
|\fA_{\fq}|\ = \
\cM_{\fq}+r(\fq)
,
\ee
say,
where
\be\label{eq:cMis}
\cM_{\fq}
:=
\sum_{\ft\mod\fq\atop\ft^{2}\equiv4(\fq)}
\sum_{\vp\in\Pi}
\frac1\fq
\sum_{q\mid\fq\atop q<Q_{0}}\sideset{}{'}\sum_{r(q)}
e_{q}(r(\tr(\vp) -  \ft))
\ee
will be treated as a ``main'' term, 
the remainder $r(\fq)$ being an error.
The analysis of the two terms is handled separately in the next two sections.


\section{Main Term Analysis}\label{sec:main}

First we wish to record the following elementary calculation. Recall that $\nu(n)$ is the number of distinct prime factors $n$.
\begin{lem}
For $\fq$ square-free,  
\be\label{eq:numFts}
\#\{\ft\in\Z/\fq\Z\ : \ \ft^{2}\equiv 4(\fq)\}
\ = \
2^{\nu(\fq)-\bo_{\{2\mid\fq\}}}
.
\ee
\end{lem}
\pf
Since $\fq$ is square-free, the equation is multiplicative. If $\fq$ is prime, then $\ft^{2}\equiv4$ implies $\ft\equiv\pm2$, 
which has two solutions unless $\fq=2$. 
\epf

We now analyze the $\cM_{\fq}$ term in \eqref{eq:cMis}, proving the following
\begin{prop}
Let $\gb$ be
the multiplicative function given at primes by
\be\label{eq:gbIs}
\gb(p) \ := \
\frac{1+\bo_{\{p\neq2\}}}p
\left(
1+{1\over p^{2}-1}
\right)
.
\ee
There is a decomposition
\be\label{eq:cMfqD}
\cM_{\fq} \ = \ \gb(\fq)\ |\Pi| + r^{(1)}(\fq)+r^{(2)}(\fq),
\ee
where 
\be\label{eq:r(1)bnd}
\sum_{\fq<\cQ}|r^{(1)}(\fq)|
\ \ll \
|\Pi|\
(\log \cQ)^{2}\
\left(
{1
\over
e^{c\sqrt{\log \faS}}}
+
Q_{0}^{C}
\faS^{-\gT}
\right)
,
\ee
and
\be\label{eq:r(2)bnd}
\sum_{\fq<\cQ}|r^{(2)}(\fq)|
\ \ll \
|\Pi|
{\cQ^{\vep}\over Q_{0}}
.
\ee
\end{prop}
\pf
Inserting the definition
\eqref{eq:PiDef} of 
$\Pi$ into \eqref{eq:cMis} gives
\beann
\cM_{\fq}
&=&
\sum_{\ft\mod\fq\atop\ft^{2}\equiv4(\fq)}
\sum_{\xi\in\Xi}
\sum_{\gw\in\gW}
\frac1\fq
\sum_{q\mid\fq\atop q<Q_{0}}\sideset{}{'}\sum_{r(q)}
\sum_{\fa\in\aleph}
e_{q}(r(\tr(\xi\fa\gw) -  \ft))
\\
&=&
\sum_{\ft\mod\fq\atop\ft^{2}\equiv4(\fq)}
\sum_{\xi\in\Xi}
\sum_{\gw\in\gW}
\frac1\fq
\sum_{q\mid\fq\atop q<Q_{0}}\sideset{}{'}\sum_{r(q)}
\sum_{\fa_{0}\in\SL_{2}(q)}
e_{q}(r(\tr(\xi\fa_{0}\gw) -  \ft))
\left[
\sum_{\fa\in\aleph\atop\fa\equiv\fa_{0}(q)}
1
\right]
.
\eeann
Apply \eqref{eq:alephED} to the innermost sum, giving
\beann
\cM_{\fq}
&=&
\cM_{\fq}^{(1)}
+
r^{(1)}(\fq),
\eeann
say, where
$$
\cM_{\fq}^{(1)}
\ := \
\sum_{\ft\mod\fq\atop\ft^{2}\equiv4(\fq)}
|\Pi|\cdot
\frac1\fq
\sum_{q\mid\fq\atop q<Q_{0}}\sideset{}{'}\sum_{r(q)}
{1\over |\SL_{2}(q)|}
\sum_{\fa_{0}\in\SL_{2}(q)}
e_{q}(r(\tr(\fa_{0}) -  \ft))
.
$$
and
\beann
|r^{(1)}(\fq)|
&
\ll
&
\tau(\fq)\
|\Pi|\
\frac1\fq
\sum_{q\mid\fq\atop q<Q_{0}}
q^{4}\
\fE(\faS;q)
.
\eeann
Here we used \eqref{eq:numFts}, and the error $\fE$ is as given in \eqref{eq:fEbnd}.
(Recall that $\tau(n)$ is the number of divisors of $n$.)
We estimate
\beann
\sum_{\fq<\cQ}|r^{(1)}(\fq)|
&\ll&
|\Pi|
\sum_{q<Q_{0}}
q^{4}\
\fE(\faS;q)
\sum_{\fq<\cQ
}
\frac{\tau(\fq)}\fq
\\
&\ll&
|\Pi|
(\log\cQ)^{2}
\bigg[
(\log \faS)^{
C
}
e^{-c\sqrt{\log \faS}}
+
Q_{0}^{C}
\faS^{-\gT}
\bigg]
,
\eeann
thus proving \eqref{eq:r(1)bnd}.

Returning to $\cM^{(1)}_{\fq}$, we add back in the large divisors $q$ of $ \fq$, writing
$$
\cM^{(1)}_{\fq}
=
\cM^{(2)}_{\fq}
+
r^{(2)}(\fq)
,
$$
say, where
$$
\cM_{\fq}^{(2)}
:=
\sum_{\ft\mod\fq\atop\ft^{2}\equiv4(\fq)}
|\Pi|\cdot
\frac1\fq
\sum_{q\mid\fq}\sideset{}{'}\sum_{r(q)}
{1\over |\SL_{2}(q)|}
\sum_{\fa_{0}\in\SL_{2}(q)}
e_{q}(r(\tr(\fa_{0}) -  \ft))
.
$$
(That is, the condition $q<Q_{0}$ has been dropped in $\cM_{\fq}^{(2)}$.)
Given $\ft$, 
let $\rho_{\ft}(q)$ be the multiplicative function given at primes by
$$
\rho_{\ft}(p):=
{1\over |\SL_{2}(p)|}
\sum_{\g\in\SL_{2}(p)}
\sideset{}{'}\sum_{r(p)}
e_{p}(r(\tr(\g)-\ft))
,
$$
so that
$$
\cM_{\fq}^{(2)}
=
\sum_{\ft\mod\fq\atop\ft^{2}\equiv4(\fq)}
|\Pi|\ 
\frac1\fq\prod_{p\mid\fq}\bigg(1+\rho_{\ft}(p)\bigg).
$$

Since $\ft^{2}\equiv4(\fq)$ and $p\mid\fq$, we have 
that
$\ft\equiv\pm2(\mod p)$.
By an elementary computation, we then evaluate explicitly that 
$$
\rho_{\ft}(p)={1\over p^{2}-1}.
$$
Using \eqref{eq:numFts}, we then have that
$$
\cM_{\fq}^{(2)}=|\Pi|\cdot \gb(\fq)
,
$$
with $\gb$ as given in \eqref{eq:gbIs}. 

Lastly, we deal with $r^{(2)}$. 
Since we crudely have
$|\rho_{\ft}(q)|\le 1/q$,  we obtain the bound
$$
|r^{(2)}(\fq)|
\ll
\tau(\fq)\
|\Pi|\
\frac1\fq
\sum_{q\mid\fq\atop q\ge Q_{0}}
\frac1q
\ll
|\Pi|\
\frac{\fq^{\vep}}\fq
\frac1
{Q_{0}}
.
$$
The estimate \eqref{eq:r(2)bnd} follows immediately, completing the proof.
\epf

\begin{rmk}
Since $\faS$ in \eqref{eq:Xis} is a small power of $N$, the first error  term  in \eqref{eq:r(1)bnd} saves an arbitrary power of $\log N$, as required in \eqref{eq:rfqLevel}. For the rest of the paper, all other error terms will be power savings. In particular, setting
\be\label{eq:Q0toga0}
Q_{0}=N^{\ga_{0}}, \qquad \ga_{0}>0,
\ee
the error in \eqref{eq:r(2)bnd} is already a power savings, while the second term in \eqref{eq:r(1)bnd} requires that
\be\label{eq:ga0togT}
\ga_{0}<{y\gT\over C}.
\ee
It is here that we crucially use the expander property for $\G$ (in fact, it is only needed for $\G_{2}\subset\G$); of course our whole point is to make the final level of distribution much larger. 
\end{rmk}


\section{Error Term Analysis}\label{sec:err}

Returning to the decomposition \eqref{eq:fAfqDecomp}, it remains to control the  error term $r(\fq)$ on average. Define
\be\label{eq:cEdef}
\cE
:=
\sum_{\fq<\cQ}|r(\fq)|
=
\sum_{\fq<\cQ}
\left|
\sum_{\ft^{2}\equiv4(\fq)}
\sum_{\pi\in\Pi}
\frac1\fq
\sum_{q\mid\fq\atop q\ge Q_{0}}\sideset{}{'}\sum_{r(q)}
e_{q}(r(\tr(\xi\fa\gw)-\ft))
\right|
.
\ee
Recall 
the decomposition $N=XYZ$ from  \eqref{eq:NXYZ}.
Our main result is the following
\begin{thm}\label{thm:cEbnd}
For any $\vep>0$, and any $1\ll 
Q_{0}<\cQ<N\to\infty$, 
\be\label{eq:cEest}
\cE\ \ll \
N^{\vep}\,
|\Pi|\,
(\xiS Z)^{1-\gd}
\left[
\frac1{Q_{0}^{1/4}}
+
{1\over Z^{1/4}}
+
{\cQ^{4}\over \xiS^{1/4}}
\right]
.
\ee
\end{thm}

As a first step, we
massage $\cE$ into a more convenient form. 
Let $\gz(\fq):=|r(\fq)|/r(\fq)$ be the complex unit corresponding to the absolute value in \eqref{eq:cEdef}, and rearrange terms as:
$$
\cE
 \ =\ 
\sum_{Q_{0}\le q<\cQ}
{1\over q}
\sideset{}{'}\sum_{r(q)}
\sum_{\vp\in\Pi}
e_{q}(r\tr(\vp))
\cdot
\gz_{1}(q,r)
,
$$
where we have set
$$
\gz_{1}(q,r) \ :=\ 
q
\sum_{\fq<\cQ\atop\fq\equiv0(q)
}
\frac{\gz(\fq)
}\fq
\sum_{\ft^{2}\equiv4(\fq)}
e_{q}(-r\ft)
.
$$
Decomposing $\Pi$ as in \eqref{eq:PiDef} and
leaving the 
special set $\aleph$ alone, we
break the $q$ sum into dyadic pieces.
This gives
\be\label{eq:cEtocE1}
\cE
\ \ll \
\sum_{\fa\in\aleph}
\sum_{Q_{0}\le Q<\cQ\atop \text{dyadic}}
\frac1Q\
|\cE_{1}(Q;\fa)|
,
\ee
where we have defined
\be\label{eq:cE1Q}
\cE_{1}(Q;\fa)
\ :=\ 
\sum_{q\asymp Q}
\left|
\sideset{}{'}\sum_{r(q)}
\gz_{1}(q,r)
\sum_{\xi\in\Xi}
\sum_{\gw\in\gW}
e_{q}(r\tr(\xi\fa\gw))
\right|
.
\ee
\thmref{thm:cEbnd} 
follows immediately
from the following estimate on $\cE_{1}(Q;\fa)$.

\begin{thm}\label{prop:cEbnd1}
We have
\be\label{eq:cE1Qbnd}
|\cE_{1}(Q;\fa)|
\ \ll  \
N^{\vep}\,
Q\,
|\Xi|
\,
|\gW|\,
(\xiS Z)^{1-\gd}
\left[
\frac1{Q^{1/4}}
+
{1\over Z^{1/4}}
+
{Q^{4}\over \xiS^{1/4}}
\right]
.
\ee
\end{thm}

\pf
To begin, 
capture the absolute value in \eqref{eq:cE1Q} by another factor $|\gz_{2}(q)|=1$,
and apply Cauchy-Schwarz in the ``long'' variable $\xi$ in \eqref{eq:cE1Q}, obtaining
$$
|\cE_{1}(Q;\fa)|^{2}
\ll 
|\Xi|
\sum_{\xi\in\SL_{2}(\Z)}
\vf_{\xiS}(\xi)
\left|
\sum_{q\asymp Q}
\gz_{2}(q)
\sideset{}{'}\sum_{r(q)}
\gz_{1}(q,r)
\sum_{\gw\in\gW}
e_{q}(r\tr(\xi\fa\gw))
\right|^{2}
.
$$
Here we have
used positivity and \eqref{eq:vfXlwr}
to
insert the weighting function $\vf_{\xiS}$ from 
\propref{prop:SL2Bnd} and
 extend the 
$\xi$ sum to all of $\SL_{2}(\Z)$. 
Since the trace of a product is a dot-product (on identifying $\Z^{4}$ with $M_{2\times2}(\Z)$ as in \eqref{eq:ZM2}), it is linear, and hence
when we open the square, we obtain
\be\label{eq:cE1Q2}
|\cE_{1}(Q;\fa)|^{2}
\ll 
\cQ^{\vep}
|\Xi|
\sum_{q,q'\asymp Q}
\sum_{\gw,\gw'}
\sideset{}{'}\sum_{r(q)\atop r'(q')}
\left|
\sum_{\xi\in\SL_{2}(\Z)}
\vf_{\xiS}(\xi)\
e\left(\xi\cdot\left[\frac rq\fa\gw-\frac{r'}{q'}\fa\gw'\right]\right)
\right|
.
\ee
Here we used crudely that
$$
|\gz_{1}(q,r)|\ll \cQ^{\vep}
.
$$

Write the bracketed expression in lowest terms as
\be\label{eq:sq0}
\frac \bs{q_{0}}:=\frac rq\fa\gw-\frac{r'}{q'}\fa\gw'
,
\ee
with $\bs=\bs(q,q',r,r',\gw,\gw',\fa)\in
\Z^{4}$ 
being coprime to
 $q_{0}\ge1$; here $q_{0}$  depends on the same parameters as $\bs$.
To study this expression in greater detail, we introduce some more notation. All variables labelled $q$, however decorated, denote square-free numbers.

Write
$$
\widetilde q:=(q,q'),\quad
q=q_{1}\widetilde q,\quad
q'=q'_{1}\widetilde q,\quad
\hat q:=[q,q']=q_{1}q_{1}'\widetilde q
,
$$
and observe from \eqref{eq:sq0} that $q_{1}q_{1}'\mid q_{0}$ and $q_{0}\mid \hat q$. Hence we can furthermore write
$$
\widetilde q_{0}:=(q_{0},\widetilde q),\quad
\hat q=q_{0}\hat q_{0}=q_{1}q'_{1}\widetilde q_{0}\hat q_{0}
,
$$
whence $q_{0}=q_{1}q_{1}'\widetilde q_{0}$. 
Note also that $Q\ll\hat q\ll Q^{2}$. 

Observe further that \eqref{eq:sq0}
implies 
$$
q_{1}'r\gw\equiv q_{1}r'\gw'\qquad\mod \hat q_{0}
,
$$
and using $\det\gw=\det\gw'=1$ gives
$$
(q_{1}'r)^{2}\equiv (q_{1}r')^{2}\qquad\mod\hat q_{0}.
$$
Since $(q_{1}r',\hat q_{0})=1=(q_{1}'r,\hat q_{0})$, we obtain
\be\label{eq:q1prhatq0}
q_{1}'r\equiv u q_{1}r'\qquad\mod\hat q_{0},
\ee
where $u^{2}\equiv 1(\hat q_{0}).$ 
There are at most $2^{\nu(\hat q_{0})}\ll N^{\vep}$ such $u(\mod \hat q_{0})$, where $\nu(m)$ is the number of distinct prime factors of $m$. 
It follows that
\be\label{eq:ggpmodq0}
\gw\equiv u\gw'\qquad\mod \hat q_{0}
.
\ee

To make full use of this last condition, we extend the $\gw$-summation to all of $\SL_{2}(\Z)$, again inserting the smoothing function $\vf$, now to parameter $Z$.
In summary, we have
\beann
|\cE_{1}(Q,\fa)|^{2}
&\ll&
|\Xi|
\sum_{Q\ll \hat q\ll Q^{2}}\
\sum_{
\begin{subarray}{l}
q_{1}q'_{1}\widetilde q_{0}\hat q_{0}=\hat q,
\\
q:=q_{1}\widetilde q_{0}\hat q_{0}\asymp Q,
\\
q':=q'_{1}\widetilde q_{0}\hat q_{0}\asymp Q,
\\
q_{0}:=q_{1}q'_{1}\widetilde q_{0}
\end{subarray}
}
\sum_{u(\hat q_{0})\atop u^{2}\equiv1(\hat q_{0})}
\sideset{}{'}\sum_{r(q)}
\sideset{}{'}\sum_{r'(q')\atop q_{1}'r\equiv uq_{1}r'(\hat q_{0})}
\sum_{
\gw'\in\gW}
\\
&&
\hskip.2in
\sum_{
\gw\in\SL_{2}(\Z),\ \gw\equiv u\gw'(\hat q_{0})
\atop
s:=q_{0}\left(\frac rq\fa\gw-\frac{r'}{q'}\fa\gw'\right)\in\cP(\Z^{4})
}
\vf_{Z}(\gw)
\left|
\sum_{\xi\in\SL_{2}(\Z)}
\vf_{\xiS}(\xi)\
e_{q_{0}}\left(\xi\cdot s\right)
\right|
.
\eeann

Working from the inside out, apply  \eqref{eq:SL2Bnd} the innermost $\xi$ sum, 
and \eqref{eq:vfXeq} to the $\gw$ sum, estimating the $\gw'$ sum trivially.
%
There are at most $q'/\hat q_{0}
$ values for $r'$, and at most $q
$ values for $r$; note that 
$$
{qq'\over \hat q_{0}}={qq'q_{0}\over \hat q}\ll  {Q^{2}q_{0}\over \hat q}.
$$
The $u$ sum contributes $N^{\vep}$, as does the sum on divisors of $\hat q$. 
Putting everything together 
gives
\beann
|\cE_{1}(Q,\fa)|^{2}
&\ll&
|\Xi|
N^{\vep}
\sum_{Q\ll \hat q\ll Q^{2}}\
\sum_{q_{0}\hat q_{0}=\hat q}
{Q^{2}q_{0}\over\hat q}
|\gW|
\bigg[
{Z^{2}\over \hat q_{0}{}^{3}}+Z^{3/2}
\bigg]
\bigg[
{X^{2}\over q_{0}^{3/2}}+q_{0}^{3}X^{3/2}
\bigg]
\\
&\ll&
N^{\vep}
Q^{2}
|\Xi|^{2}
|\gW|^{2}
(XZ)^{2(1-\gd)}
\sum_{Q\ll \hat q\ll Q^{2}}\
{1\over\hat q}
\bigg[
{1\over \hat q{\, }^{1/2}}+{1\over Z^{1/2}}
+
{
Q^{8}\over X^{1/2}}
\bigg]
,
\eeann
where we used \eqref{eq:Xisize} and \eqref{eq:gWsize}.
\thmref{prop:cEbnd1} follows immediately, as does \thmref{thm:cEbnd}.
\epf

\section{Proof of the Sieving Theorem 
}\label{sec:pfThm1}

In this section, we combine the analyses of the previous two to  prove \thmref{thm:sieve}.

Let $\fA=\{a_{N}(n)\}_{n\ge1}$ be as constructed in \eqref{eq:aNdef}. Combining \eqref{eq:fAfqDecomp} and \eqref{eq:cMfqD} gives the decomposition
$$
|\fA_{\fq}|
=
\gb(\fq)|\Pi|
+
r^{(1)}(\fq)
+
r^{(2)}(\fq)
+
r(\fq)
,
$$
as in \eqref{eq:fAfq}, with $\gb$ given by \eqref{eq:gbIs}. It is classical that \eqref{eq:gbEst} holds, so it remains to verify \eqref{eq:rfqLevel} with $\cQ$ being the level of distribution. Write
$$
X=N^{x},\ Y=N^{y},\ Z=N^{z},\ \cQ=T^{\ga}=N^{2\ga},\ Q_{0}=N^{\ga_{0}},
$$
with
\be\label{eq:xyz1}
x+y+z=1.
\ee
The bounds \eqref{eq:r(1)bnd} and \eqref{eq:r(2)bnd} are sufficient as long as $y>0$, $\ga_{0}>0$, and
\be\label{eq:yToga}
\ga_{0}\ < \ {y\gT\over C}.
\ee

The three error terms in \eqref{eq:cEest} are sufficiently controlled if
\bea
\label{eq:ga0ToX}
\ga_{0}/4&>&(x+z)(1-\gd),\\
\label{eq:zToX}
z/4&>&(x+z)(1-\gd),\text{ and }\\
\label{eq:x4Toga}
x/4&>&8\ga+(x+z)(1-\gd).
\eea
\begin{rmk}
Taking $y$ and $\ga_{0}$ very small and $x$ and $\gd$ very near $1$, it is clear that \eqref{eq:x4Toga} will not allow us to do better than $\ga<1/32$; 
 this is what we
 achieve 
 below.
\end{rmk}


Now, let $\eta>0$ be given, sufficiently small, and set
$$
\ga\ = \ \frac1{32}-\eta,
$$
as claimed in \eqref{eq:ga110}.
 We will assume further that 
$$
\gd\ > \ 1- 
\eta
$$
(more stringent restrictions on $\gd$ follow), 
and set
$$
x\ = \ 1-
\eta
,
$$
so that \eqref{eq:Xchoice} and \eqref{eq:PiSize} are satisfied.
Then 
an elementary computation shows that
\eqref{eq:x4Toga} is satisfied. 

After more elementary manipulations, we may set:
$$
z\ = \ {\eta\over 1+C/\gT}, \qquad
y \ = \ z\cdot{C\over\gT},\quad\text{ and } \quad
\ga_{0 } \ = \ \frac56 z,
$$
and assume that 
\be\label{eq:gdLwr}
\gd \ > \
1-{\eta\over 5(1+C/\gT)}
.
\ee

Then
$$
{y\gT\over C} \  =\  \frac65\ga_{0}>\ga_{0},
$$
whence \eqref{eq:yToga} is satisfied. Likewise,
$$
z/4 \ > \ \ga_{0}/4 \ = \ \frac5{24}z \ > \ \frac15z  \ > \ 1-\gd \ > \ (x+z)(1-\gd),
$$
so that \eqref{eq:ga0ToX}  and \eqref{eq:zToX} hold.
The condition \eqref{eq:gdLwr} is guaranteed by taking $\cA$ sufficiently large, cf. \eqref{eq:gdAlarge}. This completes the proof of \thmref{thm:sieve}.

\begin{rmk}\label{rmk:gTunif}
We emphasize again that it is in the last step here that we need $\aleph$ to come from the fixed group $\G_{2}\subset\G_{\cA}$. Indeed, the constants $\gT$ and $C$
are  then absolute (see \thmref{thm:BGS2}), and do not depend on $\cA$,
so \eqref{eq:gdLwr} can be ensured by taking $\cA$ large. 
\end{rmk}


\section{Proof of \thmref{thm:1}}\label{sec:pfMain}

Having established \thmref{thm:sieve} (and hence \thmref{thm:AP}) in the last section, we are in position to prove 
\thmref{thm:tr2m4}, from which we will deduce
\thmref{thm:1}.

\subsection{Proof of \thmref{thm:tr2m4}}\label{sec:pftr2}\

The deduction from \thmref{thm:AP} is straightforward, but we give the details. 
We begin by first bounding the trace multiplicity.
\begin{lem}
For any $\cA<\infty$, and any $t\ge1
$,
\be\label{eq:mulBnd}
\#\{
\g\in\G_{\cA} 
 \ : \
\tr(\g)=t
\}
\ \ll \ 
t^{1+\vep}.
\ee
\end{lem}
\pf
Let
 $\mattwos abcd\in\G_{\cA} 
 $
 have trace $a+d=t$.
 Since the entries of $\G_{\cA}$ are all positive,
 there are at most $t$ choices of $a$,
   whence $d=t-a$ is determined. Then 
   $bc=ad-1\le t^{2}$ is determined, and
   there are $\ll t^{\vep}$ choices for the divisors $b$ and $c$.
\epf

Returning to \thmref{thm:tr2m4}, let $\eta>0$ be given. Applying \thmref{thm:AP} gives a sufficiently large $\cA=\cA(\eta)$ and a set $\Pi\subset\G_{\cA}\cap B_{N}$ so that \eqref{eq:PiAPbnd} holds. 
To illustrate more clearly the mechanism below, write 
$$
\ga=1/350,
$$ 
so that
$$
\Pi_{AP}
\ = \
\{
\vp\in\Pi\ : \ p\mid (\tr(\vp)^{2}-4)
\ \Longrightarrow\  p>
N^{\ga}\}
.
$$

Now, we have
\bea
\nonumber
&&
\hskip-.5in
\#\{
\g\in\G_{\cA}\cap B_{N}\ : \ \tr(\g)^{2}-4\text{ is square-free}\}
\\
\nonumber
&
\ge
&
\#\{
\g\in\Pi_{AP} \ : \ \tr(\g)^{2}-4\text{ is square-free}\}
\\
\label{eq:nsqf}
&
>
&
N^{2\gd-\eta}
\ - \
\#
\Pi_{AP}^{\square}
,
\eea
where
we used \eqref{eq:PiAPbnd} and defined
$$
\Pi_{AP}^{\square}
\ := \
\{
\g\in\Pi_{AP} \ : \ \tr(\g)^{2}-4\text{ is not square-free}\}
.
$$
Now, 
for each $\g\in\Pi_{AP}^{\square},$ there is a prime $p$ with 
$p^{2}\mid (\tr(\g)^{2}-4)$. Since $\g\in\Pi_{AP}$, we thus have that $p>N^{\ga}$, and
moreover, 
$p^{2}$ divides either $\tr(\g)+2$ or $\tr(\g)-2$; in particular, $p\ll N^{1/2}$. 
Therefore, reversing orders and applying \eqref{eq:mulBnd}, we have
\beann
\#
\Pi_{AP}^{\square}
&
\le
&
\sum_{N^{\ga}<p\ll N^{1/2}}
\sum_{t<N\atop t^{2}-4\equiv0(p^{2})}
\#\{\g\in\G_{\cA}\cap B_{N} \ : \ \tr(\g)=t\}
\\
&
\ll
&
\sum_{N^{\ga}<p\ll N^{1/2}}
{N\over p^{2}}\
N^{1+{\vep}}
\
\ll
\
N^{2-\ga+\vep}
.
\eeann
Inserting this estimate into \eqref{eq:nsqf} gives
\beann
&&
\hskip-.5in
\#\{
\g\in\G_{\cA}\cap B_{N}\ : \ \tr(\g)^{2}-4\text{ is square-free}\}
\\
&&
>\quad
N^{2\gd-\eta}
\ - \
O(
N^{2-\ga+\vep})
.
\eeann
The above estimate is sufficient to establish \eqref{eq:tr2m4}, as long as, roughly,
\be\label{eq:gdTogA}
2\gd>2-\ga.
\ee
This explains (up to constants) the discussion on p. \pageref{eq:gaToDel} that the exponent of distribution needs to be strong enough to overcome the thinness of $\G_{\cA}$. Of course, since we have proved the above with the absolute quantity $\ga=1/350$, 
so as long as
$\gd-\eta/2>1-1/700$
 (equivalently, $\cA$ sufficiently large), we ensure that \eqref{eq:gdTogA} holds. This completes the proof of \thmref{thm:tr2m4}.
\\

We may finally present the

\subsection{Proof of \thmref{thm:1}}\

Again, this will be an easy consequence of \thmref{thm:tr2m4}. 
Let
$$
\cT
\ :=  \ 
\{t\ge1
\ : \
t^{2}-4\text{ is square-free}\},
$$
and for an integer $t
$ and $\cA<\infty$, let the trace multiplicity be
$$
\cM_{\cA
}(t)
\ := \
\#\{
\g\in\G_{\cA}
\ : \
\tr(\g)=t
\}
.
$$

Our main claim is that, for any $\eta>0$, there is an $\cA$ sufficiently large so that
\be\label{eq:numTs}
\sum_{t\in\cT\cap[1,N]}
\bo_{\{\cM_{\cA}(t)\ge t^{2\gd-1-\eta}\}}
\ >\ 
N^{2\gd-1-\eta}
.
\ee
Indeed, from \thmref{thm:tr2m4}, 
we have
\beann
N^{2\gd-\eta}
&<&
\sum_{t\in\cT\cap[1,N]}
\cM_{\cA,N}(t)
\\
&=&
\sum_{t\in\cT\cap[1,N]}
\cM_{\cA,N}(t)
\left(
\bo_{\{\cM_{\cA,N}(t)\ge W\}}
+
\bo_{\{\cM_{\cA,N}(t)< W\}}
\right)
,
\eeann
where we have introduced a parameter $W$ to be chosen later. Using \eqref{eq:mulBnd} gives
\beann
N^{2\gd-\eta}
&\ll&
N^{1+\vep}
\sum_{t\in\cT\cap[1,N]}
\bo_{\{\cM_{\cA,N}(t)\ge W\}}
\ +\ 
N\cdot
W
,
\eeann
from which \eqref{eq:numTs} follows on setting $W=N^{2\gd-1-2\eta}$, say, and renaming constants. 

Now let $\gep>0$ be given, and take $\eta>0$ small enough and $\cA$ large enough that
\be\label{eq:gepTogd}
2\gd-1-\eta > 1-\gep.
\ee
This choice of $\cA=\cA(\gep)$ corresponds to a compact region 
$$
\cY=\cY(\gep)\ \subset \
\cX\quad (\ =\ T^{1}(\PSL_{2}(\Z)\bk\bH)\ ), 
$$
as in \secref{sec:ctdf}. 
Define the set $\sD=\sD(\gep)$ to be 
$$
\sD \ := \
\{
D=t^{2}-4 \ : \
t\in\cT, \ 
\cM_{\cA}(t)>t^{2\gd-1-\eta}
\}
.
$$
All $D\in\sD$ are square-free, and hence fundamental, as are the corresponding geodesics by \lemref{lem:ELMVequiv}. Moreover,
\beann
\#(\sD\cap[1,T])
&\ge&
\#\{
t\in\cT\cap[1,\sqrt T] \ : \
\cM_{\cA}(t)>t^{2\gd-1-\eta}
\}
\\
&>&
T^{1/2-\gep}
\eeann
by \eqref{eq:numTs} and \eqref{eq:gepTogd}. This confirms \eqref{eq:sDsize}.

For each $D=t^{2}-4\in\sD$, the corresponding trace multiplicity 
satisfies
\be\label{eq:cMtBnd}
\cM_{\cA}(t)
\ >\ 
t^{1-\gep}
 \ > \ 
( \sqrt{D})^{1-\gep}
 \ \gg \ 
 |\sC_{D}|^{1-\gep}.
\ee
Now, it is {\it not} the case that each $\g\in\G_{\cA}$ corresponds uniquely to a closed geodesic on $\cX$, 
but since the corresponding visual points \eqref{eq:gaDecomp} of the geodesic are all reduced, any two differ by a cyclic permutation 
of their partial quotients. 
Recalling from \eqref{eq:ell} that the word-length metric is commensurable with the log-norm, there can be at most $C\log t$ such permutations. 
Together with \eqref{eq:cMtBnd}, this gives \eqref{eq:sCsize}, and completes the proof of \thmref{thm:1}.


\appendix

\section{Proof of 
\thmref{thm:sharp}}\label{sec:C}

In this appendix, we prove \thmref{thm:sharp}; it is a pleasure to thank Elon Lindenstrauss for suggesting the argument given here. Again, the method is more-or-less standard in the ergodic-theory community, so we only give a sketch. 

Let $\cY\subset\cX$ be a given compact region, and
let $D$ be a large non-square number (we do not require that $D$ be fundamental in this section), with corresponding class group $\sC_{D}$ and class number $h_{D}$. Recall again that Duke's theorem (now in effective form) states that, for a smooth function $\psi$ on $\cX=T^{1}(\SL_{2}(\Z)\bk\bH)$, we have that
\be\label{eq:Duke2}
\int_{\cX}\psi\ d\mu_{D}\  = \ 
\int_{\cX}\psi \ d\mu_{\cX} \ + \
O\left(D^{-c}
\cS\psi\right),
\qquad\qquad(D\to\infty),
\ee
where, as in \eqref{eq:Duke}, $\mu_{\cX}$ is Haar measure on $\cX$,  $\mu_{D}$ is the measure associated to $\sC_{D}$, namely,
$$
\mu_{D}=\frac1{h_{D}}\sum_{\g\in\sC_{D}}\mu_{\g},
$$
and $\cS\psi$ is a finite-order Sobolev norm of $\psi$ (see, e.g., \cite{ClozelUllmo2004}). The constant $c>0$ in the error rate of \eqref{eq:Duke2} could be made precise in terms of subconvexity bounds for certain $L$-functions, but we prefer to keep the exponent qualitative for ease of exposition.

Let $0\le F\le1$ be a fixed function on $\cX$ which smoothly approximates  the indicator function of the complement $\cX\setminus\cY$; in particular, we assume the support of $F$ is outside of $\cY$. Now suppose that $\g\in\sC_{D}$ is a low-lying geodesic, $\g\subset\cY$. Then, writing $T$ for the time-1 shift under the geodesic flow, we have for any $x\in\g$ that $T^{\ell}.x\in\g$, and hence
$$
\{x,T.x,\cdots,T^{k-1}.x\}\cap\supp F\ = \ \O.
$$

Let
$$
M \ := \ \int_{\cX}F\ d\mu_{\cX}
$$
be the mean of $F$, so that 
$$
F_{0}\ := \ F - M
$$
has mean zero. Furthermore, for a parameter $k$ to be chosen later (of size roughly $\log D$), define
$$
f \ := \ \frac1k\left(F_{0}+T.F_{0}+\cdots +T^{k-1}.F_{0}\right).
$$
Note that for such $x$, we have $f(x)=-M$, and hence
\beann
\frac1{h_{D}}\sum_{\g\in\sC_{D}}\bo_{\{\g\subset\cY\}}
&\le &
\mu_{D}\left(
x
 \ : \
\{x, T.x, \cdots, T^{k-1}.x\}\cap\supp F=\O
\right)
\\
&\le &
\mu_{D}\left(
x\ : \
|f(x)|\ge M
\right)
\\
&
\le
&
{1\over M^{2\ell}}
\int_{\cX}
f^{2\ell}
\
d\mu_{D}
,
\eeann
where we have introduced another parameter $1\le\ell\le k$ to be chosen later (of size a small constant times $k$). 

Apply Duke's theorem \eqref{eq:Duke2} to the last term, getting
\be\label{eq:Dapp}
\frac1{h_{D}}\sum_{\g\in\sC_{D}}\bo_{\{\g\subset\cY\}}
\ \le \
{1\over M^{2\ell}}
\int_{\cX}
f^{2\ell}
\
d\mu_{\cX}
\
+
\
O\left(
D^{-c}
C^{k}
\right)
,
\ee
where we estimated $\cS(f^{2\ell})<C^{k}$, since $F$ is fixed.
Now, the geodesic flow is a Bernouilli flow, and hence mixing of all orders. It follows that
\be\label{eq:mix}
\int_{\cX}
f^{2\ell}
\
d\mu_{\cX}
\ \ll_{F_{0}} \
\left(
{
2\ell
\over
k
}
\right)^{\ell}
.
\ee
Inserting 
\eqref{eq:mix}
into
\eqref{eq:Dapp}
gives
\be\label{eq:almost}
\frac1{h_{D}}\sum_{\g\in\sC_{D}}\bo_{\{\g\subset\cY\}}
\ \ll_{\cY} \
\left(
{
2\ell
\over
k M^{2}
}
\right)^{\ell}
\
+
\
D^{-c}
C^{k}
.
\ee
Choosing
$$
k \ = \ 
{c\over 2\log C}\cdot \log D
,
$$
say,
makes the second error in \eqref{eq:almost} of size 
$$
D^{-c}
C^{k}
\ = \ 
D^{-c/2}.
$$ 
Choosing
$$
\ell \ = \ 
{M^{2}\over 4}\cdot k,
$$
say, makes the first term in \eqref{eq:almost} of size
$$
\left(
{
2\ell
\over
k M^{2}
}
\right)^{\ell}
\ = \
\left(
{
1
\over
2
}
\right)^{\ell}
\ = \
D^{-{cM^{2}\log2/( 8\log C)}}
.
$$
This last 
 exponent determines $\gep=\gep(\cY)$,
 completing the proof.



\bibliographystyle{alpha}

\bibliography{AKbibliog}

\begin{thebibliography}{ELMV09}

\bibitem[Art24]{Artin1924}
Emil Artin.
\newblock Ein mechanisches system mit quasiergodischen bahnen.
\newblock {\em Abh. Math. Sem. Univ. Hamburg}, 3(1):170--175, 1924.

\bibitem[BGS06]{BourgainGamburdSarnak2006}
Jean Bourgain, Alex Gamburd, and Peter Sarnak.
\newblock Sieving and expanders.
\newblock {\em C. R. Math. Acad. Sci. Paris}, 343(3):155--159, 2006.

\bibitem[BGS10]{BourgainGamburdSarnak2010}
Jean Bourgain, Alex Gamburd, and Peter Sarnak.
\newblock Affine linear sieve, expanders, and sum-product.
\newblock {\em Invent. Math.}, 179(3):559--644, 2010.

\bibitem[BGS11]{BourgainGamburdSarnak2011}
J.~Bourgain, A.~Gamburd, and P.~Sarnak.
\newblock Generalization of {S}elberg's 3/16th theorem and affine sieve.
\newblock {\em Acta Math}, 207:255--290, 2011.

\bibitem[BK10]{BourgainKontorovich2010}
J.~Bourgain and A.~Kontorovich.
\newblock On representations of integers in thin subgroups of {SL}$(2,{{\bf
  {Z}}})$.
\newblock {\em GAFA}, 20(5):1144--1174, 2010.

\bibitem[BK11]{BourgainKontorovich2011}
J.~Bourgain and A.~Kontorovich.
\newblock On {Z}aremba's conjecture.
\newblock {\em Comptes Rendus Mathematique}, 349(9):493--495, 2011.

\bibitem[BK14a]{BourgainKontorovich2014}
J.~Bourgain and A.~Kontorovich.
\newblock On {Z}aremba's conjecture.
\newblock {\em Annals Math.}, 180(1):137--196, 2014.

\bibitem[BK14b]{BourgainKontorovich2014a}
Jean Bourgain and Alex Kontorovich.
\newblock On the local-global conjecture for integral {A}pollonian gaskets.
\newblock {\em Invent. Math.}, 196(3):589--650, 2014.

\bibitem[BK15]{BourgainKontorovich2015a}
Jean Bourgain and Alex Kontorovich.
\newblock The {A}ffine {S}ieve {B}eyond {E}xpansion {I}: {T}hin {H}ypotenuses.
\newblock {\em Int. Math. Res. Not. IMRN}, (19):9175--9205, 2015.

\bibitem[CU04]{ClozelUllmo2004}
L.~Clozel and E.~Ullmo.
\newblock Equidistribution des pointes de hecke.
\newblock In {\em Contribution to automorphic forms, geometry and number
  theory}, pages 193--254. Johns Hopkins Univ. Press, 2004.

\bibitem[Duk88]{Duke1988}
W.~Duke.
\newblock Hyperbolic distribution problems and half-integral weight {M}aass
  forms.
\newblock {\em Invent. Math.}, 92(1):73--90, 1988.

\bibitem[ELMV09]{ELMV2009}
Manfred Einsiedler, Elon Lindenstrauss, Philippe Michel, and Akshay Venkatesh.
\newblock Distribution of periodic torus orbits on homogeneous spaces.
\newblock {\em Duke Math. J.}, 148(1):119--174, 2009.

\bibitem[FI10]{FriedlanderIwaniecBook}
John Friedlander and Henryk Iwaniec.
\newblock {\em Opera de cribro}, volume~57 of {\em American Mathematical
  Society Colloquium Publications}.
\newblock American Mathematical Society, Providence, RI, 2010.

\bibitem[Hen89]{Hensley1989}
Doug Hensley.
\newblock The distribution of badly approximable numbers and continuants with
  bounded digits.
\newblock In {\em Th\'eorie des nombres ({Q}uebec, {PQ}, 1987)}, pages
  371--385. de Gruyter, Berlin, 1989.

\bibitem[Hen92]{Hensley1992}
Doug Hensley.
\newblock Continued fraction {C}antor sets, {H}ausdorff dimension, and
  functional analysis.
\newblock {\em J. Number Theory}, 40(3):336--358, 1992.

\bibitem[HM06]{HarcosMichel2006}
Gergely Harcos and Philippe Michel.
\newblock The subconvexity problem for {R}ankin-{S}elberg {$L$}-functions and
  equidistribution of {H}eegner points. {II}.
\newblock {\em Invent. Math.}, 163(3):581--655, 2006.

\bibitem[Hum16]{Humbert1916}
G.~Humbert.
\newblock Sur les fractions continues ordinaires et les formes quadratiques
  binaires ind{\'e}finies.
\newblock {\em Journal de math{\'e}matiques pures et appliqu{\'e}es 7e
  s{\'e}rie}, 2:104--154, 1916.

\bibitem[KS03]{KimSarnak2003}
H.~Kim and P.~Sarnak.
\newblock Refined estimates towards the {R}amanujan and {S}elberg conjectures.
\newblock {\em J. Amer. Math. Soc.}, 16(1):175--181, 2003.

\bibitem[Pop06]{Popa2006}
A.~Popa.
\newblock Central values of {R}ankin ${L}$-series over real quadratic fields.
\newblock {\em Compos. Math.}, 142:811--866, 2006.

\bibitem[Sar07]{Sarnak2007}
Peter Sarnak.
\newblock Reciprocal geodesics.
\newblock In {\em Analytic number theory}, volume~7 of {\em Clay Math. Proc.},
  pages 217--237. Amer. Math. Soc., Providence, RI, 2007.

\bibitem[Sel65]{Selberg1965}
A.~Selberg.
\newblock On the estimation of {F}ourier coefficients of modular forms.
\newblock {\em Proc. of Symposia in Pure Math.}, VII:1--15, 1965.

\bibitem[Ser85]{Series1985}
Caroline Series.
\newblock The modular surface and continued fractions.
\newblock {\em J. London Math. Soc. (2)}, 31(1):69--80, 1985.

\bibitem[SGS11]{SalehiSarnak2011}
A.~Salehi~Golsefidy and P.~Sarnak.
\newblock Affine sieve, 2011.
\newblock To appear, {\it JAMS}.

\bibitem[Zag82]{Zagier1982}
Don Zagier.
\newblock On the number of {M}arkoff numbers below a given bound.
\newblock {\em Math. Comp.}, 39(160):709--723, 1982.

\end{thebibliography}

\end{document}